\documentclass[11pt
]{article}

\usepackage{amsmath,amssymb,amstext}

\usepackage{mathbbol,mathrsfs}

\newcommand{\nek}  {\newcommand}

\nek{\vyk}[1]{}

\makeatletter
\if@twoside%
\else\fi
\makeatother

\nek{\snos}[1] {\,\footnote{~#1}}
\nek{\snom} {\,\footnotemark}
\nek{\snostext}[1] {\footnotetext{~#1}}

\nek{\ubf}{\fontseries{b}\selectfont}

\nek{\skr}{\mathscr}

\nek{\axfnt}{\fontfamily{cmss}\selectfont} 
\nek{\thfnt}{\fontseries{b}\selectfont} 

\nek{\fras}[2] {\text{\mtho\large$\TS\frac{#1}{#2}$}}
\nek{\fral}[2] {\text{\large$\frac{#1}{#2}$}}

\nek{\seci}  {\subsection}%

\newlength{\dxii}
\newlength{\dxh}
\nek{\atlh} {\addtolength{\itemsep}{\dxh}}

\nek{\itla}{\item\label}

\nek{\itsep}{\itemsep=0.6ex plus 0.2ex minus 0.15ex}
\nek{\tenu}[1]{
\itsep}

\newcounter{enuf}
\nek{\fenu}{

\addtocounter{enuf}{1}
}

\nek{\squarE}{\square}

\nek{\spnt}{\newtheorem}

\spnt{theorem}    {Theorem}[subsection]

\spnt{corollar}   [theorem]{Corollary} 
\spnt{lemm}       [theorem]{Lemma} 
\spnt{propositio} [theorem]{Proposition} 
\spnt{bag}        [theorem]{Convention} 

\spnt{numdef}     [theorem]{Definition} 
\spnt{exa}        [theorem]{Example} 
\spnt{proble}     [theorem]{Question} 
\spnt{remar}      [theorem]{Remark}  

\newtheorem{proof}{Proof}

\nek{\bbag}{\begin{bag}\rm}
\nek{\ebag}{\qed\end{bag}}
\nek{\eBag}{\end{bag}}
\nek{\bcor}{\begin{corollar}}
\nek{\ecor}{\end{corollar}}
\nek{\bdf} {\begin{numdef}\rm}
\nek{\edf} {\qed\end{numdef}}
\nek{\eDf} {\end{numdef}}
\nek{\ble} {\begin{lemm}}
\nek{\ele} {\end{lemm}}
\nek{\bexa}{\begin{exa}\rm}
\nek{\eexa}{\qed\end{exa}}
\nek{\eeXa}{\end{exa}}
\nek{\bpro}{\begin{propositio}}
\nek{\epro}{\end{propositio}} 
\nek{\brem}{\begin{remar}\rm}
\nek{\erem}{\qed\end{remar}}
\nek{\eRem}{\end{remar}}
\nek{\bte} {\begin{theorem}}
\nek{\ete} {\end{theorem}}
\nek{\bqus} {\begin{proble}\rm}
\nek{\equs} {\qed\end{proble}}

\nek{\bpF} {\begin{proof}\rm}  
\nek{\epF} {\end{proof}}

\nek{\epG} [1]{\qeG{#1}\end{proof}}

\nek{\qeD} {\hfill{${\mtho\squarE}$}}

\nek{\qed}{\qeD}

\nek{\qeG} [1]
{\hspace{0.1ex}\hfill{\text{\qed~({\sl#1\hspace{0.2ex}})}}}%
\nek{\qedD}{\hfill{$\msur\dashv\msur$}}
\nek{\envur}[1] {\begin{equation}\itsep#1\end{equation}} 
\nek{\ben}{\begin{enumerate}\itsep} 
\nek{\een}{\end{enumerate}}
\nek{\bay}{\begin{array}}
\nek{\eay}{\end{array}}

\nek{\ie} {\hbox{{\sl i.}\slsP {\sl e.}}}
\nek{\fun} {well-found\-ed}
\nek{\eg} {\hbox{{\sl e.}\slsP {\sl g.}}}
\nek{\wrt}{\hbox{w.\rmsP r.\rmsP t.}}
\nek{\noo}{\hbox{w.\rmsP l.\rmsP o.\rmsP g.}}

\nek{\slsP}{\hspace{0.15ex}}
\nek{\rmsP}{\hspace*{0.35ex}}

\nek{\axif}{\sfsl}
\nek{\axce}{\hspace{0\mathsurround}}
\nek{\axbx}[1]{\axce\protect\nolinebreak{\axfnt#1\/}%
\axce\protect\nolinebreak}
\nek{\hxbx}[1]{\axce\protect\nolinebreak{{\axfnt#1\/}\axce}%
\protect\nolinebreak}
\nek{\itea} [1] 
{\item[{\fontseries{m}\fontshape{sl}\selectfont#1\axce}%
{\fontshape{n}\selectfont:}]}

\nek{\PS}  {\axbx{Power Set}}
\nek{\Cho} {\hxbx{Choice}}
\nek{\Reg} {\hxbx{Regularity}}
\nek{\Rei} {\axbx{Regularity over\/ $\dsur\dI\dsur$}}
\nek{\SSC} {\axbx{Standard Size Choice}}
\nek{\Sep} {\hxbx{Separation}}
\nek{\Rep} {\hxbx{Replacement}}
\nek{\Sat} {\axbx{Saturation}}  
\nek{\Exn} {\hxbx{Extension}}

\nek{\tf}[1]{\text{\thfnt#1}}
\nek{\AST} {\tf{AST\/}}
\nek{\ZFC} {\tf{ZFC\/}}
\nek{\HST} {\tf{HST\/}}
 
\nek{\ttbox}[1] {{\tt{#1}}}

\nek{\cof}  {\mathop{\ttbox{cof}}}
\nek{\coi}  {\mathop{\ttbox{coi}}}
\nek{\dom}  {\mathop{\ttbox{dom}}}
\nek{\ran}  {\mathop{\ttbox{ran}}}
\nek{\card} {\mathop{\ttbox{card}}}
\nek{\Card} {\ttbox{Card}}
\nek{\st}   {\mathop{\ttbox{st}}}
\nek{\Ord}  {\ttbox{Ord}}
\nek{\galm}  {\mathop{\ttbox{gal}}} 
\nek{\pro}  [1] {\psi(#1)} 
\nek{\tmin} {\mathop{\ttbox{min}}} 
\nek{\tmax} {\mathop{\ttbox{max}}}
\nek{\tsup} {\mathop{\ttbox{sup}}} 
\nek{\prO} {\psi}

\nek{\eqr} {equivalence relation}
\nek{\Eqr} {Equivalence relation}
\nek{\aO} {{\mathord{\upa{\hspace*{0.2ex}}\Ord}}}
\nek{\aC} {{\mathord{\upa{\hspace*{0.2ex}}\Card}}}
\nek{\fin}  {^{\rbox{fin}}}
\nek{\pfin} {\pfi}
\nek{\al}{\alpha}
\nek{\da}{\delta}
\nek{\Da}{\Delta}
\nek{\ba}{\beta}
\nek{\ga}{\gamma}
\nek{\Ga}{\Gamma}
\nek{\ka} {\kappa}
\nek{\kpa} {\kappa}
\nek{\la} {\lambda}
\nek{\sg} {\sigma}
\nek{\Sg} {\Sigma}
\nek{\vpi}{\varphi}
\nek{\vt} {\vartheta}
\nek{\om} {\omega}
\nek{\omi}{\om_1}
\nek{\Om} {\Omega}
\nek{\vep}{\varepsilon}
\nek{\za} {\zeta}
\nek{\xip}{{\xi+1}}
\nek{\etp}{{\eta+1}}
\nek{\nup}{{\nu+1}}
\nek{\AL}{\aleph}
\nek{\alo}{{\AL_0}}
\nek{\ali}{{\AL_1}}
\nek{\Dzna} {{\tt ss}}
\nek{\DD}    {{\DS{{\bf\Delta}^\Dzna_2}}}
\nek{\sDD}   {{\DS{{\bf\Delta}^\Dzna_2}}}
\nek{\DDi}   {{\DS{{\bf\Delta}^\Dzna_1}}}
\nek{\DSi}   {{\DS{{\bf\Sigma}^\Dzna_1}}}
\nek{\DPi}   {{\DS{{\bf\Pi}^\Dzna_1}}}

\nek{\bbb}{\hspace{0.0ex}}
\nek{\dvoj}[1]{{\bbb{\mathbb #1}\bbb}}
\nek{\dE}{{\dvoj E}} 
\nek{\dH}{{\dvoj H}}
\nek{\dI}{{\dvoj I}}
\nek{\dL}{{\dvoj L}}
\nek{\dQ}{{\dvoj Q}}
\nek{\dR}{{\dvoj R}}
\nek{\dS}{{\dvoj S}}
\nek{\dN}{{\dvoj N}}
\nek{\dV}{{\bbb\mathbb{W}\mathbb F\bbb}}

\nek{\dlsp} {\hspace*{0.5pt}}
\nek{\dli}  {\dL[\dlsp\dI\dlsp]}

\nek{\adsp}{\hspace*{0.1ex}}
\nek{\adN} {{\mathord{\ZS^\ast\hspace{-0.30ex}\dN}}}
\nek{\adR} {{\mathord{\ZS^\ast\hspace{-0.25ex}\dR}}}
\nek{\adQ} {{\mathord{\ZS^\ast\hspace{-0.35ex}\dQ}}}
\nek{\adV} {{\mathord{\ZS^\ast\hspace{-0.25ex}\dV}}}

\nek{\aOrd} {{\mathord{\ZS^\ast\hspace{-0.25ex}\Ord}}}
\nek{\aCard} {{\mathord{\ZS^\ast\hspace{-0.25ex}\Card}}}

\nek{\skrsp}{} 
\nek{\scri}[1]{\skr #1}

\nek{\cF} {{\scri F}}
\nek{\cP} {{\scri P}}
\nek{\cX} {\scri X}
\nek{\cY} {{\scri Y}}

\nek{\pfi} {\cP_{\tt fin}}
\nek{\sta} 
{\mathord{\hspace{-0.0ex}\vphantom{|}^\SYMBS\hspace{-0.3ex}}}
\nek{\SYMBS}{{\ast}}

\nek{\ZS}{{\vphantom{X^x}}}
\nek{\sX}{{\mathord{\ZS^\ast\hspace{-0.55ex}X}}}
\nek{\sY}{{\mathord{\ZS^\ast\hspace{-0.2ex}Y}}}
\nek{\sD}{{\mathord{\ZS^\ast\hspace{-0.55ex}D}}}

\nek{\sx} {{\mathord{^\ast\hspace{-0.5ex}x}}}
\nek{\sy}   {{\mathord{^\ast\hspace{-0.35ex}y}}}
\nek{\sn} {{^\ast\hspace{-0.4ex}n}}
\nek{\sr} {{^\ast\hspace{-0.3ex}r}}
\nek{\sw} {{\mathord{^\ast\hspace{-0.4ex}w}}}
\nek{\xd} {{\mathord{^\ast\hspace{-0.5ex}d}}}
\nek{\stf}{{\ZS^\ast\hspace{-0.8ex}f}}

\nek{\sups}[2]
{\mathord{\kern 0.05em\vphantom{X}^{#2}\kern -0.17em #1}}
\nek{\ost}[1] {{\sups{#1}\circ}}
\nek{\upst} {^{\qtt{st}}}   
\nek{\qtt}[1] {\text{{\ttbox{#1}}}}

\nek{\DS}{\displaystyle}
\nek{\TS} {\textstyle}

\nek{\Askip}{\hspace{0.5pt}}
\nek{\est}{{\DS\exists\upst}}
\nek{\fst}{{\DS\forall\Askip\upst}}

\nek{\sneq}{\subsetneqq}
\nek{\sq}  {\subseteq}
\nek{\pu}  {\varnothing}
\nek{\lra} {\longrightarrow}
\nek{\ra}  {\rightarrow}
\nek{\imp} {\Longrightarrow}
\nek{\mpi} {\Longleftarrow}
\nek{\eqv} {\Longleftrightarrow}
\nek{\iy}  {\infty}
\nek{\limp}{\,\imp\,}
\nek{\leqv}{\,\eqv\,}
\nek{\res} {\mathbin{\restriction}}
\nek{\dop} [1] {{#1}\vphantom{X}^\complement}
\nek{\we}  {{\mathbin{\kern 1.3pt ^\wedge}}}
\nek{\ti}  {\times}
\nek{\dm}  {$$}
\nek{\kaz} {\forall\,}
\nek{\sus} {\exists\,}
\nek{\mto} {\mapsto}
\nek{\lmto}{\longmapsto}
\nek{\nin} {\not\in}
\nek{\bez} {\smallsetminus}
\nek{\ima} {\text{\rm"}\hspace{0.05ex}}
\nek{\imj} [2] {{#1}\hspace{0.4ex}\text{\krm''}\hspace{-0.0ex}{#2}}
\nek{\otn} [3] {{#1\mathrel{#2}#3}}
\nek{\cle} {\preccurlyeq}

\nek{\onto}{\stackrel{\rm onto}\lra}

\nek{\mem}{\dd\in}
\nek{\ste} {{\hspace{0.5pt}{\tt st}-$\mtho\in$-}}

\nek{\ang} [1]{\langle #1\rangle}
\nek{\ans} [1]{\{\hspace{0.02ex}#1\hspace{0.02ex}\}}
\nek{\sis}[2] {{\ans{#1}}\vphantom{{}^x_x}_{#2}} 

\nek{\dd} [1] {{$\rsur\mtho#1\qsur$-{\hspace{0mm}}}}
\nek{\ddd}[3] {\dd{\ang{#1,#2}}#3}
\nek{\ddt}[4] {\dd{\ang{#1,#2,#3}}#4}
\nek{\ddi}{\dd\dI}
\nek{\ddv}{\dd\dV}

\nek{\hr} {Hrba\-\v cek}
\nek{\kw} {Kawa\"\i}

\nek{\dst}{descriptive set theory}

\nek{\ek} [1] {[#1]}
\nek{\eke}[1] {[#1]_{\qE}}
\nek{\ekf}[1] {[#1]_{\qF}}

\nek{\obr}{^{-1}}

\nek{\sK}{{\sta\hspace*{-0.3ex}K}}

\nek{\fp}[2]{{\bf\Pi}^{#1}_{#2}}

\nek{\hh} {H} 

\nek{\qE} {\mathbin{\sf E}}
\nek{\qF} {\mathbin{\sf F}}

\nek{\qei} [1] {\mathbin{{\sf M}_{#1}}}
\nek{\qeu} {\qei U} 
\nek{\qeup} {\qei{U'}} 
\nek{\qev} {\qei V}
\nek{\qen} {\qei \dN}

\nek{\zo} {,\linebreak[0]}
\nek{\zi} {,\linebreak[0]\,}
\nek{\zd} {,\linebreak[0]\:}
\nek{\zt} {,\linebreak[0]\;}
\nek{\zz} {\linebreak[0]}

\nek{\mtho}{\mathsurround=0mm}
 
\nek{\mthf}{\mathsurround=0.2ex}
\nek{\msur}{\hspace*{-1\mathsurround}}
\nek{\psur}{\hspace*{1\mathsurround}}
\nek{\dsur}{\hspace{-0.3\mathsurround}}
\nek{\hsur}{\hspace{-0.5\mathsurround}}
\nek{\qsur}{\hspace{0.2\mathsurround}}
\nek{\rsur}{\hspace{0.4\mathsurround}}

\nek{\noi}{\noindent}
\nek{\vom}{\vspace{1mm plus 0.3mm minus 0.3mm}}
\nek{\vtm}{\vspace{2mm plus 0.6mm minus 0.6mm}}
\nek{\rav} [1] {{\sf D}_{#1}}

\nek{\shad} {\ost}
                
\nek{\seq} [2] {{#1}[#2]}

\nek{\ens} [2] {\ans{{#1\hspace{0.5ex}{:}}\zz\hspace{0.5ex}#2}}
\nek{\rit} [1] {{\it#1\/}}

\nek{\arin} {\stackrel{\text{\krm\tt int}}{\longrightarrow}}

\nek{\srho} {{\sta\hspace*{-0.5ex}\rho}}

\nek{\cef} [1] {|#1|^{\text{\mtho$\tt eff$}}}
\nek{\cin} [1] {|#1|^{\text{\mtho$\tt int$}}}

\nek{\eqf} {\equiv^{\text{\mtho$\tt eff$}}}
\nek{\lef} {\le^{\text{\mtho$\tt eff$}}}
\nek{\mef} {<^{\text{\mtho$\tt eff$}}}

\nek{\lej} {\le_{\text{\mtho$\tt eff$}}}
\nek{\eqj} {\equiv_{\text{\mtho$\tt eff$}}}
\nek{\mej} {<_{\text{\mtho$\tt eff$}}}

\nek{\lew} {\le^\text{\mtho$\tt w$}_{\text{\mtho$\tt eff$}}}
\nek{\eqw} {\equiv^\text{\mtho$\tt w$}_{\text{\mtho$\tt eff$}}}
\nek{\mew} {<^\text{\mtho$\tt w$}_{\text{\mtho$\tt eff$}}}

\nek{\lep} {\le^+_{\text{\mtho$\tt eff$}}}
\nek{\eqp} {\equiv^+_{\text{\mtho$\tt eff$}}}
\nek{\mep} {<^+_{\text{\mtho$\tt eff$}}}

\nek{\eqi} {\equiv^{\text{\mtho$\tt int$}}}
\nek{\lei} {\le^{\text{\tt int}}}

\nek{\spv} [1] {\|#1\|^\ast}
\nek{\spn} [1] {{\|#1\|}_\ast}

\nek{\iten} [1] {\item{\hspace*{-\labelwidth}\ubf #1}}
\nek{\itei} [1] {\item{\hspace*{-\labelwidth}\it #1}}

\nek{\ssyl} [1] {\refcite{#1}}

\nek{\lap} [1] {``#1''}

\nek{\pvn} [1] {\cP^\dI(#1)}

\nek{\krm} {\fontshape{n}\selectfont}

\nek{\upk}{\uppercase}

\nek{\bnu} {\vec \nu}

\nek{\ssi} {\hbox{s.\rmsP s.{}}}

\nek{\tsp}{\vspace{8pt plus3pt minus2pt}}

\nek{\tro} {\widehat\rho}

\nek{\vid}  {\mathop{\ttbox{wid}}}

\nek{\dt} [2] {\dd{(#1,#2)}}

\nek{\eku}  [1] {[#1]_{U}}
\nek{\ekup} [1] {[#1]_{U'}}
\nek{\ekv}[1] {[#1]_{V}}

\nek{\meu}  {\adN/U} 
\nek{\meup} {\adN/{U'}} 
\nek{\mev}  {\adN/V}

\nek{\reff} {\ref}

\nek{\refcite}{\cite}
\nek{\abstracts} [1] {\begin{abstract}#1\end{abstract}}

\begin{document}

\title{Effective cardinals in the nonstandard
universe\thanks
{This project was partially supported by
DFG grant 436 RUS 17/68/05.}
}

\author{%
Vladimir Kanovei\thanks
{IITP, Moscow, 
{\tt kanovei@mccme.ru} \ and \
{\tt vkanovei@math.uni-wuppertal.de}.
Support of RFBR 03-01-00757 acknowledged. --
{\it Contact author.}
}
\and
Michael Reeken\thanks
{Dept.\ Math.,
University of Wuppertal,
{\tt reeken@math.uni-wuppertal.de}.}
}

\date{December 2005}

\maketitle

\abstracts{
We study the structure of effective cardinals in the 
nonstandard set universe of \hr\ set theory \HST. 
Some results resemble those known in descriptive set theory 
in the domain of Borel reducibility of equivalence relations.
}

\seci*{Introduction}

Nonstandard analysis as a domain in mathematics\snos
{See \ssyl{gkk,u2}   
on the precursorial history of \rit{infinitesimal analysis}.}
emerged in the beginning of 1960s when A.\,Robinson 
\ssyl{rob'}  
demonstrated that nonstandard models 
(that is, proper elementary extensions)
of the real continuum lead to a mathematically rigorous 
system including infinitesimals and infinitely large numbers.
In the course of 1960s, the model theoretic tools used by
Robinson were shown to be applicable to variety of
mathematical structures, and that such an applicability was
based on a few general properties of nonstandard extensions,
in particular, elementarity and saturation.
For instance any \dd\ali saturated elementary extension
$\adN$ of the integers $\dN$ contains an infinitely large
number.
Several nonstandard axiomatical systems were proposed,
beginning with the mid-1970s, based on those general
principles.
Unlike the model-theoretic approach, such theories as
Nelson's internal set theory \ssyl{n77}, two theories
of \ssyl{h78,h79}, bounded set theory \ssyl{k:umn},
axiomatically described nonstandard extensions of the
whole standard set universe of \ZFC\ rather than
extensions of any  particular structure.

In the mid-1990s we formulated 
\rit{\hr\ set theory} \HST\ \ssyl{hyp},
based on earlier theories in \ssyl{h78,h79}.
This theory accumulated achievements of different
nonstandard set theories and inhibited their faults.
The set universe of \HST\ is axiomatized as a von Neumann
superstructure $\dH$ over a fully saturated elementary
extension $\dI$ ($\dI$ = internal sets)
of the class $\dV$ of all \fun\ sets,
see more on this in Section~\ref{hst}.
Our monograph \ssyl{book} presents in detail the structure
of the \HST\ universe and metamathematical properties
of \HST\ and some other popular nonstandard set theories.

This paper is devoted to the structure of
cardinalities in the nonstandard set universe of \HST.
Note that \HST\ does not include the axioms of \PS, \Cho,
and \Reg.
In fact these axioms contradict \HST. 
This is why methods of study of the structure of
cardinalities known from \ZFC\ are not always applicable
in \HST. 
Nevertheless there are two rather regular families of
cardinalities
in \HST: \ddv cardinals and \ddi cardinals.
Either family behaves in \ZFC-like manner simply because
both $\dV$ and $\dI$ satisfy \ZFC.
The intersection of the two families consists of finite
cardinals.
But little is known beyond this.
Some independence results have been obtained.
For instance,
the hypothesis that all infinite sets in $\dI$
are equinumerous in the whole universe $\dH,$
and the hypothesis that \ddi cardinals are preserved
in $\dH$
(except for hyperfinite cardinalities $m<n$ such
that $\frac mn$ is not infinitesimal, \ssyl{kkml})
are consistent with \HST, see \ssyl{ip}
or \ssyl{book}, Chapter 7.

Yet an alternative approach seems to be much more
promising in the context of \HST. 
Instead of abstract \lap{cantorial} cardinalities,
we consider here those induced by 
\rit{effective} embeddings,
\ie\ those definable in some way
or given by a certain construction.
\vyk{
This approach follows the ideology of descriptive
set theory, where similar problems are extensively
studied since the early 1990s 
(see \eg\ \ssyl{h:orb,h,ndir}).
 
In the context of nonstandard analysis,
}%
In this we follow earlier works in nonstandard analysis.
For instance studies on collapse of hyperfinite
cardinalities  by Borel and countably determined maps 
were carried out in 1980s, see
\ssyl{ast,kkml,svan}.
Further studies revealed a complicated structure of  
\lap{Borel} and
\lap{countably determined} cardinalities of hyperfinite 
sets \ssyl{mon}.

However \HST\ admits a much more general concept of effective 
cardinality than those based on Borel or countably determined 
maps. 
This concept involves the class $\dli$ of all sets 
\rit{constructible over\/ $\dI$},
and the class $\DD$ of all sets $x\in\dli\zt x\sq\dI$
(see details below), which includes and greatly exceeds 
Borel and countably determined sets.
 
The first part of the paper is devoted to effective cardinalities 
of internal sets and, generally, sets that consist of internal 
elements. 
We prove that effective cardinalities of internal sets are just 
their \ddi cardinals in the \ddi infinite domail, and resemble 
multiplicative galaxies in the hyperfinite domain. 
Effective cardinalities of $\DSi$ sets 
(\dd\dV size unions of  internal sets) 
are still linearly ordered and 
admit characterization in terms of cuts (initial segments) 
in the class $\aCard$ of all \ddi cardinals. 
Some results for cardinalities in more complicated classes 
$\DPi$ and $\DD$  will be presented, too. 

The second part of the paper considers effective cardinalities in 
their generality. 
Fortunately there is a reduction down to $\dI$: any set  in 
$\dli$ admits an effective bijection onto the quotient structure 
of the form $X/{\qE},$ where $\qE$ is a $\DD$ relation on a $\DD$ 
set $X$ (by necesity $X\sq\dI$). 

And this brings us to an analogy with modern descriptive
set theory, where cardinality problems for Borel quotient 
structures in Polish spaces  
became the focal point since early 1990s ---
especially in the form of {\it Borel reducibility\/}
of quotients and the corresponding \eqr s,
see \eg\ \ssyl{h:orb,h,ndir}.
We pursue essentially the same idea, with $\DD$ reduction
maps in the same role as Borel reductions in \dst.

Inspired by this analogy, we prove several results related to 
dichotomy of \lap{large}--\lap{small} sets, a nonstandard form 
of the Ramsey theorem, a theorem saying that quotients with rather 
small (for instance countable) classes are \lap{smooth} in a sense 
similar to the smoothness for quotients in descriptive set theory, 
and finally consider effective reducibility within the family of 
\rit{monadic} equivalence relations. 
Those readers with an experience in \dst\ may be 
interested to recognize similarities and differences with 
the set-up they are accustomed to.

\seci{Structure of the nonstandard universe}
\label{hst}

The language of \rit{\hr\ set theory} \HST\ contains two
basic predicates, the membership $\in$ and the
standardness $\st,$
hence it is called \rit{the \ste language}.
The axioms of \HST\ describe a set universe $\dH$
where the following classes are defined,
\dm
\bay{rclcl}
\dS &=& \ens{x}{\st x} &\ -&
\text{ standard sets};\\[0.8ex]

\dI &=& \ens{y}{\est x\:(y\in x)}
&\ -&
\text{ internal set};\\[0.8ex]

\dV &&  &\ -&
\text{ \fun\snom\ sets}; 
\eay
\dm%
so that $\dS\sq\dI,$ $\dI$ is an elementary extension of 
\addtocounter{footnote}{-1}%
\snostext{$\est$ and $\fst$ are shorthands for
\lap{there is a standard}, \lap{for all standard}.}%
\addtocounter{footnote}{1}%
\snostext{A set $x$ is \fun\ iff its transitive closure
has no infinite \mem decreasing chains.}%
$\dS$ in the \mem language,
$\dS$ (and $\dI$ as well) satisfies \ZFC\ in the
\mem language, the class $\dI$ is transitive, and the
universe $\dH$ is a von Neumann superstructure over $\dI.$
The universe $\dH$ satisfies all \ZFC\ axioms except for 
\Reg\ (weakened to \Rei),
\Cho\ (weakened to \SSC)
and \PS\ axioms.
The axioms of \Sep\ and \Rep\ are accepted in the \ste language.

Metamathematically, \HST\ is equiconsistent with \ZFC,
and \HST\ is a conservative extension of \ZFC\
in the sense that any \mem formula $\Phi$ is a theorem of
\ZFC\ iff $\Phi\upst$
(the relativization of $\Phi$ to $\dS$)
is a theorem of \HST.
See \ssyl{book} on axioms, metamathematics, basic
set theoretic structures, and the structure of hyperreals
in the \HST\ universe.

\bbag
\label{bl}
We argue in \HST\ below unless otherwise stated.
\ebag

{\ubf Asterisks.}
An \mem isomorphism $x\mto\sx$ of $\dV$ onto $\dS$ is
defined in \HST\ so that 
$\sx\cap\dS=\ens{\sy}{y\in x}$ for all $x\in\dV.$
The map $\ast$ is an elementary embedding of $\dV$ in $\dI$
in the \mem language. 
The classes $\dS$ and $\dV$ are \mem isomorphic and satisfy 
\ZFC.
Each of them can be unformally identified 
with the conventional set theoretic universe.
The class $\dV$ is somewhat more convenient in this role as
it is transitive and contains all its subsets,
hence some important set theoretic operations are absolute 
for $\dV$ in \HST.

{\ubf Integers and reals.\/}
The sets $\dN,\dQ,\dR$
(integers, rationals, reals)
belong to $\dV$ and are equal to resp.\ 
$(\dN)^\dV$ (\ie\ $\dN$ defined in $\dV$),
$(\dQ)^\dV,$ $(\dR)^\dV.$  
In addition $\sn=n$ for all $n\in\dN,$ therefore
$\dN\sq\adN,$
moreover $\dN$ is an initial segment in $\adN.$
The set $\adN$ coincides with the set $(\dN)^\dI$ of
all \rit{\ddi natural numbers}, 
similarly $\adQ$ and $\adR$ are equal to, resp., 
$(\dQ)^\dI$ and $(\dR)^\dI.$
Elements of $\adN,\adQ,\adR$ are often called resp.\
\rit{hyperintegers}, \rit{hyperrationals}, \rit{hyperreals}.

A hyperreal $x\in\adR$ is \rit{infinitesimal},
$x\simeq 0$ in symbols, if
$|x|<\sr$ in $\adR$ for all $r\in\dR\zt r>0,$ and 
\rit{infinitely large}, if $x\obr\simeq0,$ \ie\ 
$|x|>\sr$ for all $r\in\dR.$
A hyperreal $x$ is \rit{limited}, if
it is not infinitely large.
In this case there exists a unique $r\in\dR$ such that
$x\simeq \sr$ (that is, $x-\sr\simeq0$).
Such a real $r$ is denoted by $\shad x$
(the \rit{shadow}, or standard part, of $x\in\adR$).

{\ubf Ordinals and cardinals.\/}
The operation $\ast$ extends to proper classes
$X\sq\dV$ by $\sX=\bigcup_{x\in\dV,\,x\sq X}\sx,$ and this
does not yield contradiction provided $X\in\dV.$
Then $\adV=\dI.$
In \HST, the classes $\Card$ and $\Ord$
(all cardinals, resp., ordinals)
satisfy $\Card\sq\Ord\sq\dV$ and $\Ord=(\Ord)^\dV$
(that is, ordinals = \ddv ordinals), 
$\Card=(\Card)^\dV.$
Thus classes $\aCard\sq\aOrd\sq\dI$ are defined 
(all \ddi cardinals, resp., \ddi ordinals).
Note that $\adN\sq\aCard$.

{\ubf Sets of standard size.}
Sets equinumerous with sets in $\dV$ are called
\rit{sets of standard size}.
Note that $\card X\in\Card$ is defined then for any
set $X$ of standard size.
In \HST, sets of standard size is the same as 
well-orderable sets, 1.3.1 in \ssyl{book}.
The axiom of \Sat\ claims that every 
\dd\cap closed set $X\sq\dI\bez\ans\pu$ of standard size
has a non-empty intersection $\bigcap X.$
The axiom of \SSC\ claims the existence of a choice function
$f$ for any set $X$ of standard size
(\ie\ $f(x)\in x$ for all $x\in X\zt x\ne\pu$).
An easy consequence is the axiom of \PS\ for sets $X$
of standard size:
$\cP(X)$ is a set of standard size for any 
such $X.$  
Finite sets are sets of standard size.
On the other hand any infinite set $X\in\dI,$
for instance any set of the form $\ans{0,1,2,\dots,h},$
where $h\in\adN\bez\dN,$
is not a set of standard size.

\seci{Classes $\DD$ and $\dli$: effective sets}
\label{eff}

Which sets should be viewed as effective in \HST\,?
Following the examples of recursive, Borel, 
constructible sets, we have to choose an initial class of
sets and a set of operations applying to the initial sets.
The sets obtained this way are considered as effective.
In nonstandard set theoretic systems, 
internal sets are usually considered
as the initial sets, because of their special
role in the construction of nonstandard universes.
(In particular $\dI$  is the von Neumann basis of the \HST\
universe of sets.)
As for the operations, let us take unions and
intersections of families of standard size.

We immediately obtain the classes 
$\DSi\zt\DPi$ of all sets of the form resp.\
$\bigcup_{a\in A}X_{a}\zt \bigcap_{a\in A}X_{a},$
where $A\in\dV$ and all sets $X_{a}$ belong to $\dI,$
or, that is the same,
of the form resp.\ $\bigcup\cX\zt\bigcap\cX,$
where $\cX\sq\dI$ is a set of standard size.
(The index $^{\text{\tt ss}}$ indicates that unions and
intersections of sets of standard size are taken.)
We further define the class $\DD$ 
of all sets that can be represented both in the form 
$\bigcup_{a\in A}\bigcap_{b\in B}X_{ab},$
where $A,B\in\dV$ and all $X_{ab}$ belong to $\dI,$ 
and in the dual form
(possibly with different sets $A,B,X_{ab}$).
Note that taking, say, three operations of union and
intersection no new sets appear according to
the following result (1.4.2, 1.4.3 in \ssyl{book}).

\bpro
\label{142}
If\/ $\cX\sq\DD$ is a set of standard size then the
sets\/
$\bigcup\cX$ and\/ $\bigcap\cX$ belong to\/ $\DD.$
In addition, any set\/ $X\sq\dI$ defined in\/ $\dI$
by a\/ \ste formula with sets in\/ $\dI$ as parameters
belongs to\/ $\DD.$\qeD
\epro

Thus $\DD$ is a rather large class of sets.\snos
{\label{cd}%
There are meaningful subclasses
within $\DD,$ namely \rit{countably determined} sets, \ie\
those of the form $X=\bigcup_{b\in B}\bigcap_{n\in b}X_n,$
where $B\sq\cP(\dN)$ and all sets $X_n$ are internal
(there are different but equivalent formulations),
and \rit{Borel} sets that belong to the closure of $\dI$
under countable operations of $\bigcup$ and $\bigcap$.
These classes are considered within the model theoretic
nonstandard analysis under the assumption of \dd\ali\Sat,
\ssyl{kkml}. 
}
Yet it consists only of those sets $X$ satisfying $X\sq\dI.$
The class $\dli$ of all sets
\rit{constructible over\/ $\dI$}
extends $\DD$ on further levels of the von
Neumann hierarchy over $\dI$.   


\bdf
\label{li}
$\dli$ consists of all sets $x$ which admit a transfinite 
construction determined by a well-founded tree $T$ with sets
in $\dI$ attached to all endpoints of $T.$
The tree $T$ itself and the map which attaches internal sets
to the endpoints of $T$ belong to $\DD.$
In every node $t$ of $T$ that is not an endpoint, the set
of all sets, attached to immediate successors of $t$ in $T$
is defined.
The final set $x$ is obtained in the root of $T.$ 
\edf

Thus sets in $\dli$ are obtained via effectively coded
(in $\DD$)
transfinite iterations of the operation of assembling of
a set from its elements.
This enables us to view sets in $\dli$ as effectively
definable. 
Conversely, any effective (unformally) set
belongs to $\dli.$
Indeed it follows from theorem \ref{554}(ii) below that 
effective constructions have to be absolute for $\dli,$
hence the results of such constructions are
necessarily sets in $\dli.$  

Identifying the unformal notion of effectivity in 
\HST\ with $\dli,$ we put 
\envur{
\label{lef}
\left.
\bay{rccl}
x\lef y\,,&\text{ iff}&\text{ there is an injection }&
f\in\dli\,\text{ of }x\,\text{ into }\,y\\[0.7ex]
x\eqf y\,,&\text{ iff}&\text{ there is a bijection }&
f\in\dli\,\text{ of }\,x \,\text{ onto }\,y
\eay
\right\}.
}
and $x\mef y$ iff $x\lef y$ but $y\not\lef x.$ 
The ordinary Cantor -- Bernstein argument proves
${{x\lef y}\land{y\lef x}}\leqv{x\eqf y}$  for any
sets $x,y\in\dli.$ 
Define $\cef x,$ \rit{the effective cardinality} 
of $x\in\dli,$ to be the \dd{\eqf}equivalence class
$\ens{y\in\dli}{x\equiv^{\text{\tt eff}} y}.$
The inequalities $\cef x\le\cef y$ and
$\cef x<\cef y$ will be understood as synonimous to
resp.\ $x\lef y$ and $x\mef y$.

\bte
\label{554}
{\krm(i)}
If\/ $x\sq\dI$ then\/ ${x\in\DD}\eqv{x\in\dli}$.\vom

{\krm(ii)}
$\dli$ is a transitive class satisfying\/ \HST\,\snos 
{\label{f554}%
In fact the least class with these properties.
A suitable version of G\"odel's definition of relative
constructibility leads to exactly the same  
class $\dli$ in \HST.
See 5.5.6 in \ssyl{book}.}
and\/ $\dV\cup\DD\sq\dli.$\vom

{\krm(iii)}
For any set\/ $A\in\dli$ there is a set\/ $X\in\dI$
and an \eqr\/ $\qE$ on\/ $X,$ ${\qE}\in \DD,$
such that\/ $A\eqf{X/{\qE}}$.
\ete
\bpF
On (i), (ii) see 5.5.4 in \ssyl{book} where the class
$\DD$ is denoted by $\dE$.

(iii)
According to 5.5.4(8) in \ssyl{book}, there exist a set
$X\in\dI$ and a map ${h\in\dli}\zt{h:X\onto A}.$
Define, for $x,y\in X,$
$x\qE y$ iff $h(x)=h(y),$ and consider the map 
$a\in A\mto f(a)=\ens{x\in X}{h(x)=a}$.
\epF

Theorem~\ref{554} allows to suitably replace
$\dli$ by $\DD$ in the context of $\cef{\cdot}.$ 
For instance we conclude from \ref{554}(i) that
\eqref{lef} is equivalent to the
following in the domain of subsets of $\dI$:\pagebreak[0]
\envur{
\label{lef'}
\left.
\bay{rccl}
x\lef y\,,&\text{iff}&
\text{ there is a $\DD$ injection }&
f:x\to y \\[0.7ex]
x\eqf y\,,&\text{iff}&
\text{ there is a $\DD$  bijection }&
f:x\onto y
\eay
\right\} \text{for }\,x,y\sq\dI.
}

We begin the study of the structure of effective
cardinalities $\cef{\cdot}$ with rather simple classes,
internal sets and sets of standard size.

\seci{Effective cardinalities of internal sets}
\label{c-int}

Generally elements of $\aCard,$ that is,
\ddi cardinals, behave like
\ZFC\ cardinals since $\dI$ is a \ZFC\ universe
(in the \mem language).
Let $\cin x\in\aCard$ denote the \ddi cardinality
of a set $x\in\dI.$ 
Obviously $\cin x=\cin y$ implies $\cef x=\cef y$
since $\dI\sq\DD.$
This implication is partially reversible according to 
Corollary~\ref{ie"} below.
To figure out the effect of non-internal
maps in the domain of internal sets, let us give
some definitions.
Define, for any $x,$ 
\envur{
\label{spek}
\left.
\bay{rcr}
\spn x=\ens{\cin y}{x\supseteq y\in\dI}
&-& \text{the \rit{interior spectrum}\, of }\,x\\[0.8ex]
\spv x=\ens{\cin y}{x\subseteq y\in\dI}
&-& \text{the \rit{exterior spectrum}\, of }\,x 
\eay
\right\}.
}
Then $\spn x$ is a \rit{cut} (initial segment)
in $\aCard$ while $\spv x$ is a proper class and
a final segment in $\aCard.$
Further, for any $\ka\in\adN$ define the cuts
\dm
\ka\dN=\ens{\la\in\adN}{\sus n\in\dN\:(\la< n\ka)}\,,
\;\,
\ka/\dN=\ens{\la\in\adN}{\kaz n\in\dN\:(\la<\ka/n)} 
\dm
in $\adN,$ and the \rit{multiplicative galaxy} 
$\galm\ka=\ka\dN\bez \ka/\dN$ of $\ka.$
Then  
$\la\in\galm \ka$ \,iff\, neither of the fractions   
$\frac\ka\la\,,\,\frac\la\ka$ is infinitesimal.
To preserve the unity of notation put
$\ka\dN=\ka,$ $\galm\ka=\ans\ka$ for any
$\ka\in\aCard\bez\adN.$

Define, for $K,L\sq\aCard,$  $K\le L$ iff
$\kaz\ka\in K\:\sus\la\in L\:(\ka\le \la).$
Accordingly, $K<L$ iff $K\le L$ but $L\not\le K.$
In particular, in two cases when one of the sets
$K,L$ is a singleton, we obtain
\envur{%
\label{k<l}
\ka\le L\text{ \ iff \ }\sus\la\in L\:(\ka\le \la),
\quad\text{ and }\quad
K<\la\text{ \ iff \ }\kaz\ka\in K\:(\ka<\la).
}%
Note that galaxies are pairwise disjoint intervals in
$\aCard$  (singletons outside of $\adN$),
thus for any two galaxies $\Ga_1,\Ga_2,$
$\Ga_1<\Ga_2$ means that 
$\ka_1<\ka_2$ for any (equivalently, for all)
$\ka_1\in\Ga_1\zt\ka_2\in\Ga_2.$

See 1.4.9 and 9.6.12 in \ssyl{book}, or \ssyl{kkml},
on the next theorem. 
In the case of \ddi infinite sets the factors $\dN$
and $h$ in \ref{ie} vanish by obvious reasons. 

\bte
\label{ie}
{\krm(i)}
Suppose that\/ $X,Y\in\dI$ and\/ $f:X\to Y$ is a\/
$\DD$ map.
Then\/ $\cin X h\in\spv{\ran f}$ for any\/
$h\in\adN\bez\dN.$
In addition,\/
{\krm(a)} 
if\/ $\ran f=Y$ then\/ $\cin Y\le\cin X\dN,$
and\/
{\krm(b)} 
if\/ $f$ is an injection then\/ $\cin X\le\cin Y\dN.$\vom

{\krm(ii)}
Suppose that\/ $X\in\dI$ is infinite.
Then\/ $\cef X=\cef{X\ti\dN},$ in particular, 
$\cef Y\le\cef{X}$ for any internal\/ $Y$ with\/
$\cin Y\le\cin X\dN$.\qed
\ete

\bcor
\label{ie"}
If\/ $x,y\in\dI$ then\/ $\cef x\le\cef y$
is equivalent to\/ $\cin x\le \cin y\dN$
provided\/ $\cin y\in\adN\bez\dN$ and to just\/
$\cin x\le\cin y$ otherwise.\qed
\ecor

Thus $\cef x=\cef y$ is equivalent to 
$\galm{\cin x}=\galm{\cin y}$
in the domain $\adN\bez\dN,$
and equivalent to just $\cin x=\cin y$
outside of the domain $\adN\bez\dN.$
In the \ddi infinite domain $\aCard\bez\adN,$ the
two characterizations coincide.

\seci{Effective cardinalities of sets of standard size}
\label{c-ss}

By definition sets of standard size,
or \rit{\ssi\ sets}, are those  
equinumerous (that is, admit a bijection onto)
with sets in $\dV.$
For any \ssi\ set $X$ define
$\card X=\card W\in\Card,$ where $W$ is a set in
$\dV$ equinumerous with $X.$

\ble
\label{lss1}
{\krm(i)}
Any \ssi\ set\/ $X\sq\dI$ is\/ $\DSi$
and\/ $\DD$.\vom

{\krm(ii)}
Any \ssi\ set\/ $W$ is equinumerous with an \ssi\ 
set\/ $X\sq\dI.$\vom

{\krm(iii)}
If\/ $X\sq\dI$ is a \ssi\ set then\/
$\aCard\bez\dN\sq\spv X$ and\/ $\spv X\sq\dN.$\vom

{\krm(iv)}
If\/ $X,Y\sq\dI$ are \ssi\ sets then\/
$\card X=\card Y$ iff\/ $\cef X=\cef Y,$ thus\/
$\cef X$ can be identified with\/ $\card X.$
\ele
\bpF
{\krm(ii)}
We may assume that $W\in\dV.$ 
Then the map $w\mto\sw$ is a bijection of $W$ onto
$X=\ens{\sw}{w\in W}$
and $X$ is a set of standard size, too.

{\krm(iii)}
To prove $\aCard\bez\dN\sq\spv X$
fix $h\in\adN\bez\dN$ and apply \Sat\ 
to the family of all sets
$C_u=\ens{c\in\dI}{u\sq c\land \cin c=h},$ where
$u\sq X$ is finite.

{\krm(iv)}
Any bijection $f$ between two sets $X,Y\sq\dI$ of standard
size is itself a set of standard size, then apply (i).
\epF

It follows that \ssi\ sets are adequately
represented among $\DSi$ sets in the context of
$\card,$ and on the other hand $\cef{\cdot}$ and
$\card$ coincide on \ssi\ sets. 
The next theorem shows that
effective cardinalities of $\DD$ sets begin
with sets of standard size, where they
coincide with \fun\ cardinals, followed by the
domain of $\DD$ sets not of standard size.
It will be demonstrated below that the structure of
effective cardinalities in the latter is connected
with $\aCard$ in certain way.

\bte
\label{ss-la}
{\krm(i)}
Infinite internal sets are not \ssi\ sets.\vom

{\krm(ii)}
Any\/ $\DD$ set\/ $X$ not of standard size contains
an infinite internal subset, that is formally\/
$\dN\sneq\spn X.$%
\vom  

{\krm(iii)}
If\/ $X$ a \ssi\ set and\/ $Y$ is a\/ $\DD$ but not
\ssi\ set then\/ $\cef X<\cef Y.$
\ete
\bpF
(i)
A simple corollary of Lemma~\ref{lss1}(iii).

(ii)
By definition $\DD$ sets are \ssi\ unions of
$\DPi$ sets.
Yet it is another rather simple corollary of \Sat\
that any infinite $\DPi$ set contains an infinite
internal subset, see 1.4.11 in \ssyl{book}.


(iii)
By (ii) some number $h\in\adN\bez\dN$ belongs to
$\spn Y.$
On the other hand $h\in\spv X$ by Lemma~\ref{lss1}(iii).
This implies $\cef X\le\cef Y.$
The inequality $\cef Y\not\le\cef X$ follows from (i).
\epF

\seci{Exteriors and interiors}
\label{covin}

It turns out that internal approximations
$\spn X\zt\spv X$ are very instrumental in the study
of effective cardinalities of $\DSi$ and partly
$\DPi$ sets $X.$ 
Now a few words on cuts (initial segments) in
$\aCard$.   

\bdf
\label{cuts}
A cut $U\sq\aCard$ is \rit{standard size (\ssi) cofinal}
resp.\ \rit{coinitial}, iff there exist a 
cardinal $\vt\in\Card,$ infinite or equal to 
$1=\ans0,$ and an increasing, resp.\ decresing sequence 
$\sis{\nu_\xi}{\xi<\vt},$ of $\nu_\xi\in\aCard$
such that
$U=\bigcup_{\xi<\vt}\ens{\ka\in\aCard}{\ka<\nu_\xi},$ 
resp., $U=\bigcap_{\xi<\vt}\ens{\ka\in\aCard}{\ka<\nu_\xi}$. 
\edf

Note that \ssi\ cofinal cuts are $\DSi$ while 
\ssi\ coinitial cuts are $\DPi.$
Internal cuts, \ie\ those of the form 
$U=\ens{\ka\in\aCard}{\ka<\nu}\zt\nu\in\aCard,$ 
belong to either of the two \lap{standard size}
categories, for take 
$\vt=1$ and $\nu_0=\nu.$
See 1.4b in \ssyl{book} on the next result:

\bpro
\label{cut}
Any\/ $\DD$ cut in $\aCard$ is  \ssi\ cofinal  
or \ssi\ coinitial. 
If a cut is both 
\ssi\ cofinal and \ssi\ coinitial then it  
is internal.\qed
\epro

Coming back to $\spn X$ and $\spv X,$ note that for
any $X$ 
the intersection $\spn X\cap\spv X$ contains at most
one element.
If $\kpa\in\spn X\cap\spv X$
then there exist internal sets $Y,Z$ with $Y\sq X\sq Z$
and $\cin Y=\cin Z=\ka.$
In this case, if $\ka\in\adN$ then $X$
itself is internal with $\cin X=\ka,$  while if $\ka$
is \ddi infinite then only 
$\cef X=\cef\ka$ holds provided $X$ is $\DD.$

\ble
\label{sskof}
{\krm(i)}
If\/ $X$ is a set in\/ $\DSi,$ resp.,\/ $\DPi$
then\/ $\spn X$ is a standard size cofinal, resp.,
standard size coinitial cut in\/ $\aCard$.\vom

{\krm(ii)}
In both cases, $\spn X\cup\spv X=\aCard$.\vom

{\krm(iii)}
In both cases, if either $\spn X$ contains a largest
element $\ka,$ or $\spv X$ contains a least element
$\ka,$  then $\ka\in\spn X\cap\spv X$.
\ele
\bpF
(i)
Consider a set $\cX\sq\dI$ of standard size.
Let $X=\bigcup\cX.$
Then by \Sat\ any internal set $Y\sq X$ is covered 
by a set of the form $\bigcup{\cX'}$ where $\cX'\sq\cX$ is
finite.
On the other hand, by 1.3.3 in \ssyl{book} the set
$\pfin(\cX)=\ens{\cX'\sq\cX}{\cX'\,\text{ is finite}}$
is still a set of standard size.

Prove (ii) for $\DSi.$
Let $X=\bigcup\cX$ be as above. 
Show that any \ddi cardinal $\kpa\nin\spn X$
belongs to $\spv X.$
Take any set $Z\in\dI$ such that $X\sq Z.$
If $\cX'\sq\cX$ is finite then by definition $\bigcup{\cX'}$
is covered by an internal set of \ddi cardinality $\kpa,$
hence the set
$P_{\cX'}=
\ens{C\in\dI}{\bigcup{\cX'}\sq C\sq Z\land\cin C\le\ka}
\in\dI$
is non-empty.
Apply \Sat\ to the family of all these sets
$P_{\cX'}$.

(iii) 
Apply \Sat.
\epF

\bexa
\label{pan}
The following example\snos
{Essentially given in \ssyl{pan},
see also \ssyl{kkml}, p.\ 1172, but with a more
complicated proof based on a rather nontrivial
combinatorial theorem in \ssyl{fra}.}
of a $\DD$ set $X$ such that $\spn X\cup\spv X\sneq\aCard$ 
employs a nontrivial ultrafilter $U\in\dV$ over $\dN.$
Let $h\in\adN\bez\dN$ and $D=\ans{1,2,\dots,h}.$
The set $P=\pvn{D}=\cP(D)\cap\dI$
of all internal sets $x\sq D$
belongs to $\dI$ and satisfies $\cin P=2^h.$ 
Then\pagebreak[0]
\dm
\textstyle
U'=\ens{x\in P}{x\cap\dN\in U}
=\bigcup_{b\in U}\bigcap_{n\in b}\ens{x\in P}{n\in x}
\dm
is an ultrafilter in $P$ and a $\DD$ set.\snos
{The set $U'$ is even countably determined.}
\rit{We claim that $\spn{U'}=2^h/\dN$.}

Let $Z'\sq P$ be an internal set.
By \Sat\
(see \eg\ 9.2.15 in \ssyl{book} or 1.6 in \ssyl{kkml}),
$Z=\ens{x\cap\dN}{x\in Z'}$ is a closed subset of $U.$
It follows that the Lebesgue measure of $Z$ in $\cP(\dN)$
(identified with $2^\dN$)
is $0.$
Then easily the Loeb measure of $Z$ in $\pvn{D}$ is $0,$
so that $\cin{Z'}\in 2^h/\dN.$
Thus $\spn{U'}\sq 2^h/\dN.$
To prove the converse note that for any $u\in U$ the
set $X=\ens{x\in P}{x\cap\dN=u}$ is a $\DPi$ subset of
$U'$ that surely satisfies $\spn X=2^h/\dN$.

It follows from $\spn{U'}=2^h/\dN$ that 
$\spv{U'}=2^h$ --- by the
symmetry of the sets $U'$ and
$P\bez U'=\ens{D\bez x}{x\in P}$ within $P.$
\eexa

One can easily transform the set $U'$ as in \ref{pan}
to a $\DD$ set $X\sq\adN$
such that $\spn X=2^h/\dN$ and $\spv X=2^h\dN.$ 
The gap $\aCard\bez{(\spv X\cup\spn X)}$ consists,
in this case,
of the whole galaxy $\galm{2^h}=2^h\dN\bez 2^h/\dN$
in $\adN.$
The next theorem shows that this is a maximal possible
gap!

\bte
\label{gg}
If\/ $X$ is\/ $\DD$ and\/ $\ka\in\aCard,$
$\ka\nin\spv X\cup\spn X,$
then\/ $\ka\in\adN$ and the difference\/
$\aCard\bez{(\spv X\cup\spn X)}$ is a subset
of\/ $\galm\ka$.

Thus if\/ $X$ is\/ $\DD$ and\/ $\adN\sq\spn X$
then\/ $\spn X\cup\spv X=\aCard$.
\ete
\bpF
By definition $X=\bigcup_{a\in A}X_a$ where 
$A\in\dV$ and every $X_a$ is a $\DPi$ set.
Take any \ddi cardinal $\ka\in\spv X\bez\spn X.$
Obviously $\bigcup_{a\in A}\spn{X_a}\sq\spn X,$
thus $\ka\in\bigcap_{a\in A}\spv{X_a}$ 
by Lemma~\ref{sskof}.
It suffices to prove that any $\la\in\aCard$ belongs
to $\spv X$ in either of the two cases:
1) $\la=\ka\nin\adN,$
2) $\la\in\adN\bez\ka\dN.$
Note that $n\ka\le\la$ holds for all $n\in\dN$
in both cases.

Using \SSC, choose, for any $a\in A,$ a set $Y_a\in\dI$
such that $X_a\sq Y_a$ and $\cin{Y_a}=\ka.$
Thus $X$ is covered by the union $\bigcup_{a\in A}Y_a.$ 
For any finite $A'\sq A,$ the finite union 
$Y_{A'}=\bigcup_{a\in A'}Y_a$ is an internal set
satisfying $\cin{Y_{A'}}\le\la$ by the above.
The same application of \Sat\ as in the proof of
Lemma~\ref{sskof} yields an internal set
$Y$ still with
$\cin Y\le\la,$ satisfying $\bigcup_{a\in A}Y_a\sq Y,$
and hence $X\sq Y$ and $\la\in\aCard$.
\epF

The following corollary belongs to the
\lap{small--large dichotomy} type.
\ref{gc3} witnesses that a given $\DD$ set is rather
large \wrt\ a given cut $U$
(has rather large internal subsets),
while \ref{gc1} and
\ref{gc2} witness that $X$ is rather small
(can be covered by rather small internal sets).
The proof is easy: if $U\sneq\spn X$ then \ref{gc3}
holds by definition, otherwise apply Theorem~\ref{gg}
and get \ref{gc1} or \ref{gc2}
(or Lemma~\ref{sskof}(ii) -- in the case of
$\DSi$ and\/ $\DPi$ sets).


\bcor
\label{gc}
If\/ $X\sq\dI$ is a\/ $\DD$ set and\/ $U\sq\aCard$ is a\/
$\DD$ cut
then at least one of the following conditions holds, and
moreover\/ \ref{gc2} can be excluded for\/
$\DSi$ and\/ $\DPi$ sets\/ $X$$:$
\ben
\tenu{{\krm(A\arabic{enumi})}}
\atlh
\itla{gc1}
for any\/ $\ka\nin U$
there is an internal set\/ $Y\supseteq X$ such that\/
$\cin Y=\ka\,;$

\itla{gc2}
there exists\/ $h\in\adN\bez\dN$ such that\/
$h/\dN\sq U\sq h\dN,$ and for any\/
$\ka\in\aCard\bez h\dN$
there exists an internal set\/ $Y\supseteq X$ such that\/
$\cin Y=\ka\,;$

\addtocounter{enumi}{-1}

\tenu{{\krm(\Alph{enumi})}}
\itla{gc3}
there exists an internal set\/ $Y\sq X$ such that\/
$\cin Y\nin U$.\qed
\een
\ecor

\seci{Effective cardinalities of $\DSi$ sets}
\label{c-s1}

One may expect that the bigger $\spn X$
(or the smaller $\spv X$)
is the bigger $\cef X$ should be.
According to the next theorem, such a connection
holds for $\DSi$ sets $X$ except those satisfying
$\spn X\sq\dN$.

Following the notation in Section~\ref{c-int}, we
define, for any $K\sq\aCard,$  a cut 
$K \dN=\ens{\la}
{\sus \ka\in K\:\sus n\in\dN\:(\la\le n\ka)}$
in $\aCard$.  

\bte
\label{s1}
If\/ $X,Y$ are\/ $\DSi$ sets and\/ $\dN\sneq \spn Y$
then\/ ${\cef X\le \cef Y}$ is equivalent to\/
$\spn X\sq\spn Y \dN,$ and also to\/
$\spn X\sq\spn Y$ if\/ $\adN\sq\spn Y$.
\ete

The case $\spn Y\sq\dN$ will be considered below.

\bpF
Suppose that $X=\bigcup\cX$ and $Y=\bigcup\cY,$ where
$\cX,\cY\sq\dI$ are sets of standard size.
There is a set $D\in\dI$ such that $X\cup Y\sq D.$
Assume \noo\ that $\cX,\cY$ are \dd\cap closed families.
By \Sat, the sets of
\ddi cardinals $\ens{\cin{X'}}{X'\in\cX},$  
$\ens{\cin{Y'}}{Y'\in\cY}$
are cofinal in resp.\ $\spn X,\,\spn Y.$ 

\rit{Direction\/ $\imp$.}
Suppose otherwise.
Then there is an internal set $X'\sq X$ such that
$\cin{Y'}m<\cin{X'}/n$ for any internal $Y'\sq Y$
and $k,n\in\dN.$
As $\spn{Y}$ is a \ssi\ cofinal cut in $\aCard$
by Lemma~\ref{sskof}(i), there exists, by \Sat,
$\ka\in\aCard$ such that
$\cin{Y'}m<\ka<\cin{X'}/n$ for any internal $Y'\sq Y$
and $k,n\in\dN.$
Thus $\spn Y\dN<\ka,$ hence $\ka\in\spv Y$ by
Lemma~\ref{sskof}(ii).
In other
words, there is an internal set $Z$ such that
$Y\sq Z$ and $\cin Z=\ka.$
On the other hand, $\ka\dN<\cin{X'}$ while by
${\cef X\le \cef Y}$ there exists a $\DD$ injection
$X'\to Z,$ a contradiction to Theorem~\ref{ie}(i).

\rit{Direction\/ $\mpi$, in a stronger assumption 
that simply $\spn X\sq\spn Y.$}

\rit{Case 1}:
$\spn Y$ contains a maximal element 
$\kpa=\cin {Y_{0}},$ where $Y_0\in\cY,$
hence $Y_{0}\sq Y.$
Then for any $X'\in\cX$ the set
\dm
H_{X'}=\ens{h\in\dI}
{{h:D\to D} \,\land\,
{{h\res {X'}} \text{ is an injection }} \land\,
{\imj h{X'}\sq Y_{0}}}
\dm
is non-empty.
In addition, $H_{X''\cup X'}= H_{X''}\cap H_{X'}.$
\Sat\ yields an element $h\in\bigcap_{X'\in\cX}H_{X'}.$
Clearly $h\res X$ is an injection of $X$ into $Y_{0},$
and hence $\cef X\le \cef Y,$ as required.

\rit{Case 2}:
$\spn Y$ does not contain a maximal element, and for every
$\al\in\spn Y\cap\adN$ there exists $\ga\in\spn Y\cap\adN$
such that $\al\dN<\ga$ -- meaning that
$\ga>\al n$ for any $n\in\dN.$  
By \SSC\ there is a map 
$f:\cX\to\cY$ such that 
$\cin{X'}<\cin{f(X')}$ for all $X'\in \cX,$
and even $\cin{X'}\dN<\cin{f(X')}$  provided
$\cin{f(X')}$ (then also $\cin{X'}$) belongs to $\adN.$
Then  
\dm
H_{X'}=\ens{h\in\dI}
{{h:D\to D} \,\land\,
{{h\res {X'}} \text{ is an injection }} \,\land
{\imj h{X'}\sq f(X')}}
\dm
is non-empty for any  $X'\in\cX.$ 
Then argue as in  Case 1.

\rit{Case 3}:
the negation of cases 1, 2.
Then there is a number $c\in\adN\bez\dN$ such that 
$c\in\spn Y$ but $2c\nin\spn Y.$
Then $\cef{[0,c)}\le\cef Y$ while 
$\cef X\le \cef{[0,2c)}$ (see case 1).
However $\cef{[0,2c)}=\cef{[0,c)}$  by
Corollary~\ref{ie"}.

\rit{Direction\/ $\mpi$, general case\/}.
If $\spn X\sq\spn Y\dN,$
but $\spn X\sq\spn Y$ does {\ubf not} hold
then there exist numbers $c\in\adN\bez\dN$ and $n\in\dN$
such that $\spn X\sq{[0,nc)}$ and
${[0,c)}\sq \spn Y\sq{[0,2c)}.$
We have $\cef X\le \cef{[0,nc)}$ by the above, and
$\cef{[0,c)}\le \cef{Y}.$
It remains to apply Corollary~\ref{ie"}.
\epF

It remains to consider the case $\spn Y\sq\dN$ avoided
in the theorem.
It leads to sets of standard size!

\ble
\label{p.ss}
For a set\/ $X\sq\dI$ to be of standard size each of
the conditions\/ $\spn X\sq\dN,$ $\adN\bez\dN\sq\spv X$
is necessary and, if\/ $X$ is\/ $\DD$, also
sufficient.
\ele
\bpF
By Theorem~\ref{ss-la}(i)
$\spn X\sq\dN.$
On the other hand $\adN\bez\dN\sq\spv X$ by
Lemma~\ref{lss1}(iii).
The sufficiency follows from Theorem~\ref{ss-la}(ii).
\epF

Thus Theorem~\ref{s1} \rit{fails} in the case
$\spn X=\dN$:
take any pair of infinite sets $X,Y\sq\dI$ of
standard size
with $\card X\ne\card Y$ and apply 
\ref{p.ss} to show that
$\spv X=\spv Y=\dN,$ and Lemma~\ref{lss1} to show that
$\cef X\ne\cef Y.$
Nevertheless we easily obtain the following corollary.

\bcor
\label{s1c}
If\/ $X,Y$ are\/ $\DSi$ sets then their effective
cardinalities are comparable in the sense that
at least one of the following inequalities holds$:$
$\cef X\le \cef Y$ or\/ $\cef Y\le \cef X.$\qed
\ecor

\seci{Effective cardinalities of $\DPi$ sets}
\label{cp1}

The proof of $\imp$ in Theorem~\ref{s1} does not work
for $\DPi$ sets since $\spn Y$ is now \ssi\ coinitial
and the \Sat\ argument does not work.
On the other hand there is a suitable counterexample.

\bexa
\label{-pi}
Fix $h\in\adN\bez\dN$ and
let $S$ be the set of all internal maps
$s:\ans{0,1,2,\dots, h-1,h}\to\ans{0,1}=2.$
Define $a_s,b_s\in 2^\dN$ (hence $\in\dV$) so that
$a_s(k)=s(k)$ and $b_s(k)=s(h-k)$ for all $k\in\dN.$
For $a,b\in2^\dN$ put $S_{ab}=\ens{s}{a_s=a\land b_s=b}$
and $S_a=\ens{s}{a_s=a}.$
Then $S$ is internal, $\cin S=2^{h+1},$ while each 
$S_a$ is a $\DPi$ set with $\spn{S_a}=2^h/\dN.$
Obviously $(2^h/\dN)\dN= (2^h/\dN).$ 
To see that $S$ and $S_a$ lead to a counterexample to
$\imp$ of Theorem~\ref{s1}, it suffices to prove that
$\cef S=\cef{S_a}$ for some $a.$

Since either of $S,S_a$ is a union
of $2^\dN$-many sets of the form $S_{ab},$ 
it remains to show that $\cef{S_{ab}}=\cef{S_{a'b'}}$
for all $a,b,a',b'.$
By \Sat\ there is $\sg\in S$ such that
$a(n)=a'(n)\oplus \sg(n)$ and
$b(n)=b'(n)\oplus\sg(h-n)$
for all $n\in\dN,$ where $\oplus$ is addition
modulo 2.
Finally the internal map $s\mto s\oplus\sg$
(in the termwise sense)
easily maps $S_{ab}$ onto $S_{a'b'}$ in 1-1 way.
\eexa

In fact \ref{-pi} is the only possible counterexamle
for $\DPi$ sets in the following sense:
if $X,Y$ are $\DPi$ sets, $\dN\sneq\spn Y,$ and
$\cef X\le\cef Y$ then either $\spn X\sq\spn Y\dN$ or
there is a number $\ka\in\spn Y\zt\ka\in\adN\bez\dN,$
such that $\spn X=\ka/\dN$ while $\spn Y\sq\ka\dN.$
We skip the proof.

Our further goal is to present what looks like a
near-counterexample, (ii) of Theorem~\ref{p1},
to $\mpi$ of Theorem~\ref{s1} in the field of
$\DPi$ sets.

If $X$ is a $\DPi$ set then $\spv X$ is standard size
coinitial by Lemma~\ref{sskof}.
If $\spv X$ contains a least element $\ka$ then $\ka$
is simultaneously the largest element in $\spn X$
still by Lemma~\ref{sskof}, and then easily
$\cef X=\cef\ka.$
It follows that if in this case $Y$ is another
$\DPi$ set with $\spv Y=\spv X$ then $\cef X=\cef Y.$
But if $\spv X$ does not contain a least element then
there is an infinite coinitial sequence with standard
size many terms. 
This case is considered by the next theorem.
It follows from (ii) that there are sets of the
largest effective cardinality among all $\DPi$ sets
$X$ with the same $\spv X,$ while (iii) presents a
rather nontrivial partial counterexample to
Theorem~\ref{s1} for $\DPi$ sets.
We deal with \ddi infinite cardinals here, but similar
results can be obtained in the hyperfinite domain
--- we leave it to the reader.


\bte
\label{p1}
{\krm(i)}
If\/ $\cX,\cY\sq\dI$ are sets of standard
size, $X=\bigcap\cX\zt Y=\bigcap\cY,$
$\vt=\card\cX\in\dV$ is an infinite regular cardinal,
$\spv X=\spv Y,$ and the coinitiality of\/ $\spv X$
is exactly\/ $\vt,$ then\/ $\cef Y\le\cef X.$\vom

{\krm(ii)}
There exist\/ $\DPi$ sets\/ $X,Y$ as in\/ {\krm(ii)} 
such that\/ $\cef X\le\cef Y$ fails via\/ $\DD$
injections\/ $g$
of the form\/
$g=\bigcup_{w\in W}\bigcap_{\xi<\vt}g_{w\xi},$  where
all\/ $g_{w\xi}$ are internal and\/ $W$ is a set
of standard size.
\ete
\bpF
(i)
Assume \noo\  that there exist sets
$X_0\in\cX\zt Y_0\in\cY$ such that
$X\sq X_0$ and $Y\sq Y_0$
for all $X\in\cX\zt Y\in\cY,$
and the families $\cX,\cY$ are \dd\cap closed. 
We claim that there exists 
a function $\psi:\cX\to\cY$ satisfying
\envur{
\label{psi}
\kaz A\in\pfin(\cX)\:\sus f\in F\:\kaz X\in A\:
(\psi(X)\sq\dom f\land \imj f{\psi(X)}\sq X),
}
where $F\in\dI$ is the set of all 1--1 functions 
$f\in\dI$ with $\dom f\sq Y_0$ and $\ran f\sq X_0.$
To define $\psi$ fix an enumeration
$\cX=\ens{X_\al}{\al<\vt}.$ 
Suppose that $\al'<\vt,$ and the values
$\psi(X_\al)\in\cY\zt\al<\al',$ have been defined.  
In our assumptions, there is a set $Y\in \cY$
such that $\cin Y< \cin{\bigcap_{\al\in A}X_\al}$
for every finite $A\sq [0,\al'].$
To complete the inductive step put $\psi(X_{\al'})=Y.$

To prove \eqref{psi} consider a finite set
$A=\ans{\al_1<\dots<\al_n}\sq\vt.$
By the construction
$\cin{\psi(X_{\al_k})}<\cin{\bigcap_{1\le i\le k}X_{\al_k}}$
for all $k=1,\dots,n.$
Arguing in $\dI,$ we easily find a map $f\in F$ such
that
$\imj f{(\psi(X_{\al_k}))}\sq
\bigcap_{1\le i\le k}X_{\al_k}$
for every $k=1,\dots,n,$ hence \eqref{psi} holds.

Yet by \Sat\ \eqref{psi} is equivalent to the following:
\envur{
\label{psi'}
\sus f\in F\:\kaz X\in\cX\,
(\psi(X)\sq\dom f\land\imj f{(\psi(X))}\sq X).
}
Thus $\imj f{Y}\sq X,$ for such an $f,$ and hence 
$\cef Y\le\cef X$ holds even by means of an internal
map $f$.

(ii)
Fix an infinite cardinal $\vt$ in $\dV.$
It easily follows from \Sat\ that there exists a strictly
decreasing sequence $\bnu=\sis{\nu_\xi}{\xi<\vt}$ of
\ddi cardinals $\nu_\xi\in\aCard\bez\adN.$
A \rit{\dd\bnu large set} will be any $X\in\dI$
such that $\sus\xi<\vt\:(\cin X\ge\nu_\xi).$ 

Let $\tau=\vt^+$ (the next cardinal in $\dV$).
The counterexample wis based on 
a sequence $\sis{Y_\ga}{\ga<\tau}$ of
internal sets $Y_\ga$ such that
\ben
\tenu{(\alph{enumi})}
\itla{x1}
for any pair of disjoint finite sets\/
$u,v\sq\tau\zt u\ne\pu,$
the set
$Y_{uv}=\bigcap_{\al\in u}Y_\al\bez\bigcup_{\ba\in v}Y_\ba$
is\/ \dd\bnu large$\,;$

\itla{x2}
$\cin{Y_{uv}}=\cin{Y_{u\pu}}$ for any disjoint finite\/
$u,v\sq\tau\,;$

\itla{x3}
if\/ $\xi<\vt\zt A\sq\tau,$ and\/
$\cin{Y_{\ans{\al,\ba},\pu}}\ge\nu_\xi$
{\krm(that is, $\cin{Y_\al\cap Y_\ba}\ge\nu_\xi$)}
for all\/
$\al,\ba\in A$ then\/ $\card A\le\vt$.
\een

We define $Y_\ga$ by induction.
To begin with put $Y_0=[0,\nu_0)$
(an initial segment in $\aOrd$).
Now suppose that $\ga<\tau$ and a set $Y_\da\in\dI$
has been defined for every $\da<\ga$ so that
\ref{x1} and \ref{x2} hold below $\ga.$ 
Re-enumerate
$\ens{Y_\da}{\da<\ga}=\ens{Z_\al}{\al<\la},$
where $\la=\tmin\ans{\ga,\vt},$ without repetitions.

For any pair of disjoint finite sets
$u,v\sq\la\zt u\ne\pu,$ define the internal set
$Z_{uv}=\bigcap_{\al\in u}Z_\al\bez\bigcup_{\ba\in v}Z_\ba.$
In our assumptions, the \ddi cardinals
$\ka_{uv}=\cin{Z_{uv}}$ satisfy $\ka_{uv}=\ka_u,$
where $\ka_u=\ka_{u\pu},$ and
$\sus\xi<\vt\,(\ka_u\ge\nu_{\xi}).$
For any finite $u\sq\la$ let $\xi(u)$ be the least
ordinal $\xi<\vt$ such that $\nu_\xi<\ka_{u}$
and $\xi>\tsup u.$
We assert that there is an internal set $Z$ satisfying
\ben
\tenu{(\alph{enumi})}
\addtocounter{enumi}3
\itla{+}
$\cin{Z\cap Z_{uv}}=\nu_{\xi(u)}$ and 
$\cin{Z_{uv}\bez Z}=\ka_{u}$
for any pair of disjoint finite sets $u,v\sq\la,$
$u\ne\pu,$ and

\addtocounter{enuf}1
\itla{++}
$\cin{Z\bez \bigcup_{\ba\in v}Z_\ba}=\cin Z\ge\nu_0$
for each finite set $v\sq \la$.
\een
Indeed as $\vt$ is a set of standard size it suffices
to prove that for any finite $d\sq\vt$ there is a
set $Z\in\dI$ satisfying \ref+, \ref{++} for all
$u,v\sq d$.

Note that the sets of the form $Z_{uv},$ where
$u\cup v=d$ and $u\cap v=\pu,$ are mutually disjoint,
and by definition satisfy
$\nu_{\xi(u)}\le \ka_{u}=\cin{Z_{uv}}.$
This allows us to define an internal $Z$ satisfying
\ref+ for all pairs $u,v$
with $u\cup v=d\zd u\cap v=\pu\zt u\ne\pu,$ and, adding a
sufficient portion out of $\bigcup_{\ba\in d}Z_\ba,$
also $\cin{Z\bez \bigcup_{\ba\in d}Z_\ba}=\cin Z\ge\nu_0.$
It remains to show \ref+ for all disjoint
sets $u,v\sq d$ not necessarily with $u\cup v=d.$

We show this by backward induction on the
cardinality of $u\cup v.$ 
Suppose that $u\cup v\sneq d.$
Take any $\al\in d\bez{(u\cup v)}.$
Let $u'=u\cup\ans\al$ and $v'=v\cup\ans\al.$
Then by the inductive hypothesis
$\cin{Z\cap Z_{u'v}}=\nu_{\xi(u')}$ and
$\cin{Z\cap Z_{uv'}}=\nu_{\xi(u)}.$
Since $Z_{uv}= Z_{u'v}\cup Z_{uv'}$ and
easily $\xi(u)\le\xi(u')$ whenever $u\sq u',$
we conclude that 
$\cin{Z\cap Z_{uv}}=\nu_{\xi(u)}$ as required.
Similarly, $\cin{Z_{u'v}\bez Z}=\ka_{u'}$ and
$\cin{Z_{uv'}\bez Z}=\ka_{u},$ therefore
$\cin{Z_{uv}\bez Z}=\ka_{u'}+\ka_u=\ka_u$ as
required.

Take as $Y_\ga$ any set $Z\in\dI$ satisfying
\ref+, \ref{++}.
We have to demonstrate that \ref{x1}, \ref{x2} remain
true for the sequence
$\sis{Y_\da}{\da\le\ga},$ or, that is equivalent,
for the sequence $\sis{Z_\al}{\al\le\la},$
where $Z_\la=Y_\ga=Z$.

Take any pair of disjoint sets $u,v\sq\la\cup\ans\la.$
If $\la\nin u\cup v$ then the set
$Z_{uv}=\bigcap_{\al\in u}Z_\al\bez\bigcup_{\ba\in v}Z_\ba$
is the same as above so there is nothing to prove. 
Suppose that $\la\in u;$ put $u'=u\bez\ans\la.$     
Then $Z_{uv}=Z\cap Z_{u'v},$ and hence $Z_{uv}$ is
\dd\bnu large by \ref+ (applied for the pair $u',v$).
Separately if $u=\ans\la$ then $u'=\pu,$ hence
$Z_{u'v}$ is not defined, but obviously
$Z_{uv}=Z\bez \bigcup_{\ba\in v}Z_{\ba},$
therefore $Z_{uv}$ is \dd\bnu large by \ref{++}.
Suppose that $\la\in v;$ put $v'=v\bez\ans\la.$
Then $Z_{uv}=Z_{uv'}\bez Z,$ and hence $Z_{uv}$ is
\dd\bnu large still by \ref+.
This proves \ref{x1}; the derivation of \ref{x2}
from \ref+, \ref{++} is similar.

This ends the recursive construction of the sets $Y_\ga$.

Show that such a sequence $\sis{Y_\ga}{\ga<\tau}$
also satisfies \ref{x3}.
We prove not only that $\card A\le\vt$ for any set
$A$ as in \ref{x3}, but even more the order type of
$A$ in $\tau$ is $\le\vt.$ 
Suppose that $\ga\in A.$
Let us come back to the reenumerated system 
$\ens{Y_\da}{\da<\ga}=\ens{Z_\al}{\al<\la},$
where $\la=\tmin\ans{\ga,\vt},$ and to the
construction of $Y_\ga=Z_\la=Z$ satisfying \ref+.
It follows from \ref+ that, for any $\al<\la,$
$\cin{Y_\ga\cap Z_{\al}}=\nu_{\xi(\ans\al)}<\nu_\al.$
In other words, for any $\xi<\vt$ the inequality
$\cin{Y_\ga\cap Z_{\al}}\ge\nu_\xi$ can be true only
for $\al<\xi.$
Thus there exist \dd{({<}\vt)}many sets $Y_\da\zt\da<\ga,$
satisfying $\cin{Y_\ga\cap Y_{\da}}\ge\nu_\xi,$ as
required.

Coming back to the proof of (ii) of Theorem~\ref{p1},
we fix a sequence $\sis{Y_\ga}{\ga<\tau}$ satisfying
\ref{x1}, \ref{x2}, \ref{x3}, 
and put $\cY=\ens{Y_\ga}{\ga<\tau}.$
Then $Y=\bigcap_{\ga<\tau}Y_\ga$ is a $\DPi$ set.
Note that every \ddi cardinal $\nu_\xi\zt\xi<\vt,$
belongs to $\spv Y$ by \ref{x3}, and on the other
hand it follows by \Sat\ that every
internal superset $H$ of $Y$ contains a subset of the
form $\bigcap_{\al\in u}Y_\al=Y_{u\pu},$ where
$u\sq\tau$ is finite, and hence $\cin H\ge\nu_\xi$
for some $\xi<\vt$ by \ref{x1}.
It follows that the sequence
$\sis{\nu_\xi}{\xi<\vt}$ is coinitial in $\spv Y.$
It follows from Lemma~\ref{sskof} that $\spn Y$
coincides with the set
$\Om=\ens{\ka\in\aCard}{\kaz\xi<\vt\:(\ka<\nu_\xi)}$.

The other side of the counterexample will be the
$\DPi$ set $X=\bigcap_{\xi<\vt}X_\xi,$ where
$X_\xi=\ens{\ka\in\aOrd}{\ka<\nu_\xi}\in\dI.$
Easily $\cin{X_\xi}=\nu_\xi,$ therefore
the sequence $\sis{\nu_\xi}{\xi<\vt}$ is coinitial
in $\spv Y,$ too.
We conclude that $\spv X=\spv Y,$ hence
$\spn X=\Om=\spn Y$ by Lemma~\ref{sskof}.

To accomplish (ii),
suppose towards the contrary that there is
an injection $g\in\DD\zt g:X\to Y$ of the
form $g=\bigcup_{w\in W}\bigcap_{\xi<\vt}g_{w\xi},$
where all\/ $g_{w\xi}$ are internal and 
$W$ a set of standard size.
Then each $g_w=\bigcap_{\xi<\vt}g_{w\xi}$ is still
an injection into $Y,$ whose domain
$D_w=\dom g_w\sq X$ is still a $\DPi$ set, moreover,
an intersection of \dd{({\le}\vt)}many internal sets.
(The combination of quantifiers
$\sus\:\fst\xi<\vt$ converts to
$\fst p\in\pfin(\vt)\:\sus$ by \Sat.)

We claim that $\spn{D_w}=\Om$ for at least one
$w\in W.$

(Indeed otherwise choose any
$\ka_\al\in\Om\bez\spn{D_w}$ for every $w\in W;$
here \SSC\ is applied.
Recall that $\spv X$ is a standard size coinitial
final segment in $\aCard,$ therefore the complement
$\Om$ of  is not standard size cofinal by \Sat.
It follows that there is an \ddi cardinal $\ka\in\Om$
bigger than each $\ka_w.$
Then $\ka\in\spv{D_w}$ for any $w\in W,$ thus
any $D_w$ is covered by an internal
set of \ddi cardinality $\ka.$
Still by \Sat, the union $\bigcup_{w\in W}D_w$
can be covered by an internal set $C\zt \cin C=\ka.$
Then $X\sq C,$ contradiction.)

This result allows us to replace $X$ by $D_a,$ or,
in different words, reduce the task to the case when
$g,$ a given injection $X\to Y,$ is equal to
$\bigcap_{\xi<\vt}g_\xi,$ each $g_\xi$ being an
internal set.
An easy application of \Sat\ shows that there  is a
finite set $u\sq\vt$ such that
$h=\bigcap_{\xi\in u}g_\xi$ is an injective function.
On the other hand $h$ is an internal function extending
$g,$ thus $X=\dom g\sq h$ and $\imj hX\sq Y.$
Let $D=\dom h$ (an internal superset of $X$).

Recall that $Y=\bigcap_{\ga<\tau}Y_\ga$ where all
$Y_\ga$ are internal.
Thus, for any $\ga,$ $\imj hX\sq Y_\ga,$ and hence, as
$X=\bigcap_{\xi<\vt}X_\xi$ and the family of all sets
$X_\xi$ is \dd\cap closed, \Sat\ yields an ordinal
$\xi(\ga)<\vt$ such that $X_{\xi(\ga)}\sq D$ and
$\imj h{X_{\xi(\ga)}}\sq Y_\ga.$
As $\tau=\vt^+,$
there is at least one $\xi<\vt$ such
that $G=\ens{\ga<\tau}{\imj h{X_{\xi}}\sq Y_\ga}$
is unbounded in $\tau.$
Thus all sets $Y_\ga\zt\ga\in G$ include as a subset
one and the same internal set $R=\imj h{X_{\xi}}.$
Note that $\cin{R}=\cin{X_{\xi}}$ because $h$ is an
injection.
But $X_\xi$ is a \dd\bnu large set, a contradiction
with \ref{x3}.                                 
\epF

\bqus
\label{mh}
Is Corollary~\ref{s1c} still
true for sets in $\DD$ or in  $\DPi$\,?
\equs

Theorem~\ref T below shows that a wider category of
$\DD$ \rit{quotients} has plenty of incomparable sets.
Note that the existence of countably determined
sets incomparable in the sense of countably determined
injections, is also an open problem.
A counterexample defined in \ssyl{cv} in the \AST\
frameworks makes use of the hypothesis
that there exist only \dd\ali many internal sets, and
hence is irreproducible in \HST.

On the other hand all Borel sets
(in the sense of Footnote~\reff{cd})
are Borel-comparable.
This result was first obtained by \AST-followers,
see \eg\ \ssyl{ast}, and then reproved in \ssyl{svan}.
See more on this in \ssyl{book}, 9.6 and 9.7.

\seci{Effective sets in the form of quotients}
\label{qg}

Sets of the form $X/{\qE},$ where $X$ is $\DD$ while
$\qE$ is a $\DD$ equivalence relation on $X$ will be
called \rit{$\DD$ quotients}.
These $\DD$ quotients include the class
$\DD$ itself, for take $\qE$ to be just the equality
on a given $\DD$ set $X,$ so that the map sending any
$x\in X$ to $\ans x$ is
a bijection of $X$ onto ${X/{\qE}}.$

On the other hand, it follows from Theorem~\ref{554}(iii)
that every set in $\dli,$ that is, every effective set in
the sense explained in Section~\ref{eff}, admits an
effective bijection onto a $\DD$ quotient.
Thus $\DD$ quotients exhaust, in the context of effective
cardinalities, all effective (= $\dli$) sets in general.

One may ask whether $\DD$ quotients produce more
effective cardinalities than just $\DD$ sets.
Call \rit{smooth} any $\DD$ quotient that admits a
$\DD$ bijection onto a $\DD$ set.
We show in Section~\ref{ctble}
that every $\DD$ quotient ${X/{\qE}},$ such that all
\dd\qE classes $\eke x=\ens{y\in X}{x\qE y}\zt x\in X,$
are sets of standard size, is smooth.
A family of non-smooth $\DD$ quotients, those defined
by means of \rit{monadic partitions} of $\adN,$ will
be studied in Sections~\ref{mop}, \ref{redu}.
We prove there that there exist incomparable effective
cardinalities of monadic $\DD$ quotients, still an
open problem for $\DD$ sets themselves. 
We also prove a \lap{small--large} type theorem for
$\DD$ quotients in Section~\ref{fm},
similar to \ref{gc} but not so sharp, with an interesting
Ramsey-like corollary.

Note that $\DD$ quotients consist of subsets of
$\dI$ which are not necessarily internal sets themselves.
Accordingly injections of $\DD$ quotients are maps
whose $\dom$ and $\ran$ not necessarily consist of
internal sets.
Still there is a way to pull the consideration down
to the basic level. 

\bdf
\label E
Let ${\qE}\zi{\qF}$ be \eqr s\ on sets $X,Y.$
A set $R\sq X\ti Y$ is a 
\rit{\dt\qE\qF invariant pre-injection of\/ $X$ into\/
$Y$} iff
\,1)
$\dom R=X$\snos
{This condition can be weakened to
$\eke x\cap \dom R\ne \pu$ for any $x\in X$
without any harm.}
and
\,2)
the equivalence ${x\qE x'}\eqv {y\qF y'}$ holds
for all
$\ang{x,y}\in R$ and $\ang{x',y'}\in R.$

Such a set $R$ is a
\rit{reduction of\/ $X/{\qE}$ to\/ $Y/{\qF}$}
(or just of $\qE$ to\/ $\qF$)
if in addition
\,3)
$R$ is a (graph of a) function $X\to Y.$ 

Write $\qE\lej\qF$ iff there is a
\dt\qE\qF invariant pre-injection
$P\sq X\ti Y,$ $P\in\DD,$ of $X$ into $Y.$
Write $\qE\lep\qF,$ in words: 
\rit{$\qE$ is effectively reducible to\/ $\qF,$}
iff there is a reduction $\rho\in\DD,$
$\rho: X\to Y$ of $\qE$ to $\qF.$

An \eqr\ $\qE$ on a set $X$ and the quotient
$X/{\qE}$ are \rit{\dd\DD smooth} iff there is a 
$\DD$ set $Y$ such that $\qE\lej\rav Y,$ where
$\rav Y$ is the equality on $Y$ considered as
an equivalence relation.~\snos
{\label{smut}%
Note that in this case any invariant pre-injection
is a partial map that can be immediately extended
to a reduction, and hence in fact $\qE\lep\rav Y$ holds.}
\edf

This definition resembles some central concepts in
modern descriptive set theory, like Borel reducibility
and \lap{Borel cardinals}
(see, for instance, \ssyl{h:orb,h,ndir}), where
Borel maps are used in approximately the same role as
$\DD$ maps in this paper.

\bpro
\label R
{\krm(i)}
Suppose that\/ ${\qE}\zi{\qF}$ are\/ $\DD$ \eqr s\ on\/
$\DD$ sets\/ $X,Y.$
Then\/ $\cef{X/{\qE}}\le \cef{Y/{\qF}}$ iff\/
$\qE\lej\qF$.

{\krm(ii)}
An $\DD$ \eqr\ $\qE$ on a\/ $\DD$ set\/ $X$ is\/
\dd\DD smooth
iff there exists a\/ $\DD$ set\/ $Y$ such that\/
$\cef{X/{\qE}}=\cef Y.$ 
\epro
\bpF
(i)
Suppose that $f\in\dli$ is an injection
${X/{\qE}}\to{Y/{\qF}}.$
Then
$P=\ens{\ang{x,y}\in X\ti Y}{f(\eke x)=\ekf y}$
is a set in $\dli,$ hence a $\DD$ set by \ref{554}(i),
and obviously an invariant pre-injection. 
The converse is equally simple:
if $P$ is an invariant pre-injection then to define
an injection $f:{X/{\qE}}\to{Y/{\qF}}$ put 
$f(\eke x)=\ekf{y}$ for any $\ang{x,y}\in P.$

(ii)
Suppose that $\qE\lej\rav Z,$ where $Z$ is a $\DD$ set.
Let this be witnessed by an invariant pre-injection
$R\sq X\ti Z$ of class $\DD.$
Clearly $R=\rho$ is then a reduction
(a map $X\to Z$ such that
${x\qE x'}\eqv {\rho(x)=\rho(x')}$).
The set $Y=\ran \rho\sq Z$ is as required.
\epF

\seci{\Eqr s with standard size classes}
\label{ctble}

In modern descriptive set theory, an \eqr\ $\qE$  is 
\rit{countable} iff all equivalence classes 
$\eke x=\ans{y:x\qE y}\zt x\in\dom\qE,$
are at most countable.
See \ssyl{ctble} on properties and some open problems
related to countable \eqr s.
But in the nonstandard setting the structure of \eqr s in
a much wider class turns our to be considerably simpler:
all of them admit effective transversals.

Recall that a
\rit{transversal} of an \eqr\ is any set having exactly
one element in common in every equivalence class.      

\bte
\label{sq}
Any\/ $\DD$ \eqr\/ $\qE,$ on an internal set\/ $\hh$
and with \ssi\ classes, has a\/ $\DD$ transversal
and hence is\/ \dd\DD smooth. 
\ete

Opposed to this, the Vitali equivalence on the reals is
obviously countable but not smooth (via Borel maps), neither
it admits a Borel transversal.


\bpF 
First of all, a $\DD$ transversal implies \dd\DD smoothness:
let $\rho(x)$ denote the only element of the transversal
equivalent to $x$ and apply \ref{142} to show that $\rho$
is still $\DD.$
Let us prove the existence of a $\DD$ transversal.

By definition
${\qE}=\bigcup_{a\in A}\bigcap_{b\in B}E_{ab},$
where $E_{ab}\sq\hh\ti\hh$ are internal sets while
$A,B\in\dV.$
Put $\imj Px=\ens{y}{\ang{x,y}\in P}$ for  
$P\sq\hh\ti\hh$ and $x\in\hh.$

\ble
\label{sq1}
There exists a standard size family\/ $\cF$
of internal maps\/ $F:\hh\to\hh$ such that\/
$\eke x\sq\ens{F(x)}{F\in\cF}$ for all\/ $x\in\hh.$
\ele
\bpF
It suffices to prove the lemma for each \lap{constituent}
$E_a=\bigcap_{b\in B}E_{ab}$ of $\qE. $
According to 1.3.6 in \ssyl{book}, 
the intersection $\bigcap\cX$ of a \ssi\ family $\cX$
of internal sets either is not a \ssi\ set or it is finite
and there is a finite $\cX'\sq\cX$ such that
$\bigcap\cX'=\bigcap\cX.$
It follows that every set $\seq{E_a}{x}$ is finite 
and moreover there is a finite set
$\ba_{ax}\sq B$ such that
$\imj{E_a}{x}=\bigcap_{b\in \ba_{ax}}\imj{E_{ab}}{x}.$
Put, for any $n\in\dN$ and any finite $\ba\sq B,$
\dm
\textstyle
E_{a\ba}=\bigcap_{b\in \ba}E_{ab}
\quad\text{and}\quad
P_{a\ba n}=\ens
{\ang{x,y}\in E_{a\ba}}
{\card\imj{E_{a\ba}}{x}\le n}\,.
\dm
All sets $P_{a\ba n}$ are internal.
We define
$F_{a\ba ni}(x)=$ \lap{\dd i\:th element of $P_{a\ba n}$ in 
the sence of a fixed internal linear ordering of
$P_{a\ba n}$} in the case when $1\le i\le n$ and
$P_{a\ba n}$ contains at least $i$ elements,
and $F_{a\ba ni}(x)=y_0$ otherwise,
where $y_0$ is a once and for
all fixed element of $\hh.$ 
It remains to define $\cF$ to be the family of all functions
$F_{a\ba ni}$.
\epF

Let $\cF$ be as in the lemma. 
The sets
\dm
D_{F}=\dom {(\qE\cap F)}=\ens{x\in\hh}{x\qE{F(x)}}
\qquad
(F\in\cF)
\dm 
belong to $\DD$ by Proposition \ref{142}.
Let us fix an internal wellordering $\prec$ of the set
$\hh.$
%
Suppose that $F\in\cF.$ 
For any $x\in\hh$ we carry out the following construction 
called the \rit{\dd Fconstruction for\/ $x$.} 
Define an internal \dd\prec decreasing sequence   
$\sis{x_{(a)}}{a\le a(x)}$ of length $a(x)+1\in\adN.$   
Its terms $x_{(a)}$ are defined by induction on $a.$
Put $x_{(0)}=x.$ 
Assume that $x_{(a)}$ has been defined. 
If $z=F(x_{(a)})\prec x_{(a)}$ then put $x_{(a+1)}=z,$
otherwise put $a(x)=a$ and stop the construction. 
Eventually the construction ends since   
$x_{(a+1)}\prec x_{(a)}$ for all $a.$
Put $\nu_{F}(x)=0$ if $a(x)$ is even and $\nu_{F}(x)=1$ 
otherwise.

Define $\pro x(F)=\nu_{F}(x)$ for any $x\in\hh\zt F\in\cF;$
thus $\prO:\hh\to 2^\cF.$

\ble
\label{Zr}
If\/ $r\in2^\cF$ then\/
$\Psi_r=\ens{x\in\hh}{\pro x=r}$ belongs to\/ $\DD$.
\ele
\bpF
Note that $x\in\Psi_r$ iff  
$\nu_F(x)=r(F)$ for all $F\in\cF.$
On the other hand, all sets
$X_{F}=\ens{x\in\hh}{\nu_F(x)=0}$ ($F\in\cF$) 
are internal because the \dd Fconstruction is internal.
It remains to apply Proposition~\ref{142}.
\epF

According to the next lemma, any two different but
\dd\qE equivalent elements $x\in\hh$ 
have different \lap{profiles} $\pro x$.

\ble
\label{ple'}
If\/ $x\ne y\in\hh$ and\/ $x\qE y$ then\/ $\pro x\ne \pro y$.
\ele
\bpF
Suppose that $y\prec x.$ 
There exists a function $F\in\cF$ such that $y=F(x).$ 
Then $y=x_{(1)}$ in the sense of \dd{F}construction for $x.$ 
It follows that the \dd{F}construction for $y$ has exactly 
one step less than the \dd{F}construction for $x.$
Thus $\nu_{F}(x)\ne\nu_{F}(y)$ and $\pro x\ne \pro y$.
\epF

We continue the proof of Theorem~\ref{sq}.
Note that $2^\cF$ and $\cP(2^\cF)$ are sets of standard
size together with $\cF$ (1.3.3 in \ssyl{book}).
Thus by the axiom of \SSC\ there is a map
$A\mto r_A$ such that $r_A\in A$ for any non-empty 
$A\sq 2^\cF.$
Its graph $C=\ens{\ang{A,r}}{A\sq 2^\cF\land r=r_A}$
is a \ssi\ set together with $\cP(2^\cF).$ 
For any $x\in\hh$ put $A(x)=\ans{\pro y:y\in\eke x},$
a non-empty subset of $2^\cF.$
Now
$
X=\ans{x\in\hh:\pro x=r_{A(x)}}
$
is a transversal for $\qE$ by Lemma~\ref{ple'}.

Prove that $X$ is a $\DD$ set.
By definition $X=\bigcup_{\ang{A,r}\in C}Y_{A}\cap\Psi_r,$
where
$
Y_{A}
=\ens{x\in \hh}{A(x)=A}.
$
However $\Psi_r\in\DD$ by Lemma~\ref{Zr}.
It remains to check that $Y_{A}\in\DD$
for each $A\sq 2^\cF.$ 
Note that
$A(x)=\ens{\pro{F(x)}}{{F\in\cF}\land{x\in D_F}},$ 
and hence $A(x)=A$ is equivalent to
\dm
\bay{l}
\kaz r\in A\:\sus F\in\cF\:
({x\in D_F}\land {F(x)\in \Psi_r})
\;\,\land\,\hspace*{16ex}\\[0.8ex]
\hspace*{16ex}\land\;\,\kaz F\in\cF\:\sus r\in A\:
({x\in D_F}\imp {F(x)\in \Psi_r})\,.
\eay
\dm
Yet the sets $\Psi_r$ and $D_F$ are $\DD$ 
(see above),   
while the domains $A$ and $\cF$ are sets of standard 
size.
Now apply Proposition~\ref{142}.\vom

\epG{Thm~\ref{sq}}

\seci{Monadic partitions}
\label{mop}

A cut $U\sq\adN$ is \rit{additive} if
${a\in U}\imp{2a\in U}.$
Any such cut $U$ induces an equivalence 
relation $x\qeu y$ iff $|x-y|\in U$ on $\adN.$
(The additivity implies that $\qeu$ is transitive.)                                              
Its equivalence classes 
$
[x]_U=\ans{y:x\qeu y}=\ans{y:|x-y|\in U},
$ 
are called \dd U{\it monads\/} and relations of the form
$\qeu,$ accordingly, \rit{monadic} \eqr s or
\rit{monadic partitions}.

Monads of various kinds are considered in nonstandard
analysis.
As for those induced by additive cuts in $\adN,$ see
\ssyl{jin,KL}. 

The following is an elementary corollary of
Proposition~\ref{cut}:

\bpro
\label 2
If\/ $\pu\ne U\sneq\adN$ is an additive\/ $\DD$ cut
then\/ $U$ is non-internal and either standard size
cofinal or standard size coinitial.\qed
\epro

Any additive $\DD$ cut $U\sq\adN$ defines a $\DD$
quotient $\meu=\adN/{\qeu},$ the set of all \dd Umonads.
According to the next theorem, effective cardinalities  
of those quotients are determined by two factors.
The first of them~is 
\dm
\TS
\vid U\,=\, 
\bigcap_{u\in U,\,u'\in\adN\bez U}[0\,,\fras{u'}u) 
\,=\,
\bigcap_{u\in U}\bigcup_{u'\in U,\;u'>u}
[0\,,\fras{u'}u),
\dm
\rit{the width} of $U.$\snos
{Also called {\it the thickness\/} of $U$ in some
papers on \AST.}
The second one is the cofinality/coinitiality.
The \rit{cofinality} $\cof U$   
of a standard size (\ssi) cofinal non-internal cut,
is the least cardinal $\vt\in\Card$ such that
$U$ has an increasing cofinal sequence of type $\vt.$
The \rit{coinitiality} $\coi U$ of a standard size
coinitial
cut is defined similarly, with a reference to coinitial
sequences in $\aCard\bez U.$
Note that $\cof U$ and $\coi U$ are infinite regular
cardinals.

Additive cuts of lowest possible width are obviously   
those of the form $U=c\dN\zt c\in\adN$ and 
$U=c/\dN\zt c\in\adN\bez\dN,$ which we call {\it slow\/}; 
they satisfy $\vid U=\dN.$
Other additive cuts will be called {\it fast\/}.  

\bte
\label T
Suppose that\/ $U,\,V$ are additive\/ $\DD$ cuts in\/ 
$\adN$ other than\/ $\pu$ and\/ $\adN.$ 
Then\/ {\krm(i)} $\cef{\adN}\le\cef{\adN/U}.$
In addition$,$
\ben
\tenu{{\krm(\roman{enumi})}}
\addtocounter{enumi}1
\atlh
\itla{T0}
$\adN/U$ is\/ \dd\DD smooth iff\/
$\adN/U$ has a\/ $\DD$ transversal
iff\/ $U$ is slow$;$

\itla{T1}
if\/ $U$ is slow then\/ $\cef{\adN/U}\le\cef{\adN/V}\,;$

\itla{T2f}
if both\/ $U,\,V$ are \ssi\ cofinal cuts and\/ $U$
is fast then\/ $\cef{\adN/U}\le\cef{\adN/V}$\,
iff\/$:$\,
$\cof U=\cof V$ and\/ $\vid U\sq\vid V\,;$
 
\itla{T2i}
if both\/ $U,\,V$ are \ssi\ coinitial cuts and\/
$U$ is fast then\/ $\cef{\adN/U}\le\cef{\adN/V}$\,
iff\/$:$\,
$\coi U=\coi V$ and\/ $\vid U\sq\vid V\,;$

\itla{T4}
if\/ $U,V$ are fast cuts, $U$ is \ssi\ cofinal
and\/ $V$ is \ssi\ coinitial then\/
$\cef{\adN/U}$  and\/ $\cef{\adN/V}$ are incomparable.
\een
\ete

Thus either of the two classes of monadic partitions
(\ssi\ cofinal and \ssi\ coinitial)
is linearly 
\dd\lej(pre)ordered in each subclass of the same
cofinality (coinitiality), slow partitions of both classes
form the \dd\lej least type, and there is no other
\dd\lej connection between the two classes and their
same-cofinality/coinitiality subclasess. 

See \ssyl{mon} for earlier results of countably
determined and Borel reducubility of monadic
partitions for \rit{countably} cofinal/coinitial cuts.

\seci{The proof of the reducibility theorem}
\label{redu}

We begin the proof of Theorem~\ref T with the
following observation.

\brem
\label{scar}
Call a set $X\sq\adN$ \rit{scattered} iff there is
a number $c\in\adN\bez\dN$ such that
$\fras{\cin{X\cap I}}c$ is infinitesimal for any
interval $I$ in $\adN$ of length $c.$
It is quite clear that $\adN$ is not
a finite union of scattered sets, and hence, by \Sat,
\rit{$\adN$ is not a standard size union of
internal scattered sets}.
\erem

\noi{\ubf Proof of Theorem~\ref T.}
(i)
Choose a number $h\in\adN\bez U.$
The map $x\mto \eku{xh}$ is an 
injection of $\adN$ into $\adN/U$.  

\ref{T0}
If $\meu$ admits a $\DD$ transversal
then it is \dd\DD smooth.
(Let, for $x\in\adN,$ $\rho(x)$ be the only element
of the transversal equivalent to $x.$)
Suppose that $\meu$ is smooth, \ie\
$\qeu\lej\rav R$ for a suitable $\DD$ set $Z.$
This is witnessed by a $\DD$ reduction
$\rho:\adN\to Z$
By Theorem~\ref{ie}(i) the set $\ran\rho$ can be
covered by an internal set $Y$ with
$\cin Y\le\cin\adN.$
Thus $\cef{\meu}\le\cef\adN.$
Then $\cef{\meu}\le\cef\mev$ for any other additive
$\DD$ cut $V$ by (i), thus $U$ must be slow by \ref{T4}.
Finally, if $U$ is slow then $\meu$ has a
$\DD$ transversal by Theorem~1.4.7 in
\ssyl{book}.\snos
{Theorem~\ref{sq} yields a $\DD$ transversal for
$\adN/\dN,$ and hence for any
$\adN/{(h\dN)}$ by multiplication.
Transversals defined this way are countably
determined but not Borel.
Yet partitions of the form
$\adN/{(h/\dN)}$
have no countably determined transversals
by 9.7.14 in \ssyl{book}.
These theorems were obtained in \ssyl{book}
on the base of earlier results in \ssyl{jin}.}

\ref{T1}
If $U$ is slow then $\meu$ is \dd\DD smooth, and in
fact $\cef\meu\le\cef\adN,$ see the proof of \ref{T0}.
It remains to apply (i).

\ref{T2f} 
Thus let $U\zi V$ be additive \ssi\ cofinal cuts.
Choose increasing sequences $\sis{u_\xi}{\xi<\vt}$
and $\sis{v_\eta}{\eta<\tau}$ 
cofinal in resp.\ $U$ and $V;$ $\vt=\cof U$
and $\tau=\cof V$ being
infinite regular cardinals in $\dV.$
As $U$ is supposed to be fast, we can assume that
$\fras{u_{\xip}}{u_\xi}$
is infinitely large for all $\xi.$ \vom 

\rit{Part 1}: assuming $\cef\meu\le\cef\mev,$
we prove that $\vid U\sq\vid V.$ 
Let, by \ref R,
$R\sq \adN\ti\adN$ be a \dt UVinvariant pre-injection,
thus $\dom R=\adN,$ 
and
${|x-x'|\in U}\eqv{|y-y'|\in V}$ for all pairs 
$\ang{x,y}$ and $\ang{x',y'}$ in $R.$

Since $R$ is $\DD,$ we have, by definition, 
$R=\bigcup_{a\in A}\bigcap_{b\in B}R_{ab},$ 
where $A,B\in\dV$ and the sets
$R_{ab}\sq\adN\ti\adN$ are internal.

Let us fix $a\in A.$

Then $R_a=\bigcap_{b\in B}R_{ab}\sq R,$ hence for
any $\eta<\tau$
we have\pagebreak[0] 
\dm
\kaz b\:(
{\otn x{R_{ab}}y}\,\land\, 
{\otn{x'}{R_{ab}}{y'}}
)
\,\land\,
{|y-y'|<v_\eta}\,\limp\,
\sus \xi<\vt\:(|x-x'|<u_\xi) 
\dm
for all $x,x',y,y'\in\adN.$
We obtain, by \Sat, 
\envur{
\label{eq1}
\bay[b]{l}
\hspace*{-2ex}
\kaz\eta<\tau\quad\sus\text{ finite }F\sq B\;\:
\sus \xi<\vt\quad
\kaz x,x',y,y'\in\adN:\\[0.4\dxii]
\qquad{\otn x{R_{aF}}y}\,\land\, 
{\otn{x'}{R_{aF}}{y'}}
\,\land\,{|y-y'|<v_\eta}\,\limp\,
{|x-x'|<u_{\xi}}\,,
\eay
}
where $R_{aF}=\bigcap_{b\in F}R_{ab}.$ 
A similar (symmetric) argument yields:
\envur{
\label{eq2}
\bay[b]{l}
\hspace*{-2ex}
\kaz\xi<\vt\quad\sus\text{ finite }F'\sq B\;\:
\sus \eta<\tau\quad
\kaz x,x',y,y'\in\adN:\\[0.4\dxii]
\qquad{\otn x{R_{aF'}}y}\,\land\, 
{\otn{x'}{R_{aF'}}{y'}}
\,\land\,{|x-x'|<u_\xi}\,\limp\,
{|y-y'|<v_{\eta}}\,.
\eay
}

Suppose, towards the contrary, that 
$\vid U\not\sq\vid V.$
Then there exists $\eta<\tau$ such that the sequence
$\sis{\fras{v_{\eta'}}{v_\eta}}{\eta<\eta'<\tau}$
is not cofinal in $\vid U.$

\rit{Keeping\/ $a\in A$ still fixed}, we let
$F$ and $\xi$ satisfy \eqref{eq1} for this $\eta.$  
By the choice of $\eta,$ there exists an ordinal 
$\xi'>\xi$ such that
$\fras{u_{\xi'}}{u_\xi}>\fras{v_{\eta'}}{v_\eta}$
for any $\eta'>\eta,$ hence in fact
$\fras{u_{\xi'}}{u_\xi}>\ell\cdot\fras{v_{\eta'}}{v_\eta}$ 
for any $\eta'>\eta$ and any $\ell\in\dN.$
We now let 
$F'$ and $\eta'$ satisfy \eqref{eq2} (as $F$ and $\eta$)
for the $\xi'$ considered. 
We may assume that $F\sq F'$ and
$\eta'\ge \eta$ --- otherwise take, resp., the union
and the maximum of the two. 
Then we have, for all 
$\ang{x,y}\zt\ang{x',y'}$ in the set
$R(a)=R_{aF'}:$\pagebreak[0]
\envur{
\label{eq3}
\left.
\bay{lcl}
{|y-y'|<v_\eta}\,&\,\limp\,&\,
{|x-x'|<u_{\xi}}
\\[0.7\dxii]
{|x-x'|<u_{\xi'}}\,&\,\limp\,&\,
{|y-y'|<v_{\eta'}}
\eay
\right\};\quad
\xi,\xi',\eta,\eta'\,\text{ depend on }\,a.
}
Put $D(a)=\dom {R(a)},$ an internal subset of $\adN$
together with $R(a)$.

Note that any interval of length $v_{\eta'}$ in $\adN$
consists of approximately $s=\fras{v_\eta'}{v_\eta}$
subintervals of length $v_\eta.$
Accordingly any interval of length $v_{\xi'}$  
consists of approximately $t=\frac{u_\xi'}{u_\xi}$
subintervals of length $u_\xi,$ while
$\frac st$ is infinitesimal by the above.
It follows by \eqref{eq3} that
$\fras{\cin{I\cap D(a)}}{\cin{I}}$ is infinitesimal
for any interval $I$ in $\adN$ of length $u_\xi',$
hence $D(a)$ is scattered in the sense of \ref{scar}.

On the other hand
$\adN=\dom R=\bigcup_{a\in A}D_a=\bigcup_{a\in A}D(a),$
where $D_a=\dom R_a,$  simply
because $R_a\sq R(a),$ which is a contradiction
with \ref{scar}.
\vom

\rit{Part 2}:
in the same assumptions and notation as in Part 1,
we prove that $\cof U=\cof V.$ 
This means to prove $\vt=\tau.$ 
Suppose $\vt\ne\tau.$
Let say $\vt<\tau.$
(The other case is similar.)
Then, for a fixed $a\in A,$
there is an ordinal $\eta<\tau,$ one and the
same for all $\xi<\vt,$ such that
\eqref{eq2} takes the~form:
\envur{
\label{eq2'}
\bay[b]{l}
\hspace*{-2ex}
\kaz\xi<\vt\quad\sus\text{ finite }F'\sq B\;\quad
\kaz x,x',y,y'\in\adN:\\[0.4\dxii]
\qquad{\otn x{R_{aF'}}y}\,\land\, 
{\otn{x'}{R_{aF'}}{y'}}
\,\land\,{|x-x'|<u_\xi}\,\limp\,
{|y-y'|<v_{\eta}}\,.
\eay
}
Take an ordinal $\xi<\vt$ for this $\eta$ by \eqref{eq1},
and then apply \eqref{eq2'} for $\xip.$
We obtain a finite set $F\sq B$ such that, for all 
$x,x'\in D(a)=\dom{R_{aF}}:$
\envur{
\label{eq3'}
{|x-x'|<u_{\xip}}\,\limp\,{|x-x'|<u_{\xi}}.
}
However, as $U$ is fast, the cofinal sequence
$\sis{u_\xi}{}$ can be chosen so that
$\fras{u_\xi}{u_\xip}$ is infinitesimal for all $\xi.$
Then the set $D(a)$ is scattered by \eqref{eq3'}, and so
on towards the contradiction as in Part 1.\vom

\rit{Part 3}.
Suppose that $\cof U=\cof V=\vt$
(an infinite regular cardinal in $\Card$) and
$\vid U\sq \vid V.$
To prove $\cef\meu\le\cef\mev$ it suffices, by \ref R,
to define a reduction
of $\meu$ to $\mev.$ 
Let $\sis{u_\xi}{\xi<\vt}\zt\sis{v_\xi}{\xi<\vt}$
be increasing cofinal sequences
in the cuts resp.\ $U,V.$
Due to additivity of the cuts,
we may \noo\ assume that all terms $u_\xi,v_\xi$
are powers of $2$.

We first define subsequences of the cofinal sequences
satisfying a certain term-to-term
inequality.
Note that $\vid U\sq \vid V$ basically means
\dm
\kaz v\in V\quad\sus u\in U\quad\kaz u'\in U,\,u'>u\quad
\sus v'\in V,\,v'>v\;
\left(\fras{u'}u\le\fras{v'}v\right)\,.
\dm
This allows us to define an unbounded subsequence of
$\sis{u_\xi}{\xi<\vt}$ such that, after the reenumeration,
the following holds
($\xi,\eta,\za$ are ordinals $<\vt$):\pagebreak[0]
\dm
\kaz\za\quad
\kaz\xi>\za\quad
\sus\eta>\za\;
\left(
\fras{u_{\xi}}{u_\za}\le\fras{v_{\eta}}{v_{\za}},
\quad\text{that is},\quad
\fras{v_{\za}}{u_\za}\le\fras{v_{\eta}}{u_{\xi}}
\right)
\,,
\dm
and then to once again define an unbounded
subsection of, now, $\sis{v_\eta}{\eta<\vt}$
to satisfy, after the reenumeration,
the following:\pagebreak[0]
\envur{
\label{uv}
\kaz \xi<\eta<\vt\;
\left(
\fras{v_{\xi}}{u_\xi}\le\fras{v_{\eta}}{u_{\eta}},
\quad\text{that is},\quad
\fras{u_{\eta}}{u_\xi}\le\fras{v_{\eta}}{v_{\xi}}
\right)\,.
}

Finally, we may assume that $u_0=1.$
(Replace each $u_\xi$ by
$u'_\xi=\fras{u_\xi}{u_0}.$
As all $u_\xi$ are  powers of $2,$ these fractions
belong to $\adN.$
The sequence $\sis{u'_\xi}{}$ is then cofinal in the
cut $U'=U/u_0=\ens{u}{uu_0\in U}.$
The inequality
$\cef\meu\le\cef\meup$
is witnessed by the map
$\eku x\mto
\ekup{\text{entire part of }\fras x{u_0}}$.)

Note that the map $f$ sending each $u_\xi$ to $v_\xi$
satisfies the following:
$\dom f=\ens{u_\xi}{\xi<\vt}$ is a \ssi\ set,
$\dom f$ and $\ran f$ consist of powers of $2,$
and $\fras{f(u)}{u}\le\fras{f(u')}{u'}$ for all
$u<u'$ in $\dom f$ by \eqref{uv}.
By \Sat\ there is an internal
function $F$ with $D=\dom F$ a hyperfinite subset of
$\adN\bez\ans0,$ such that $\dom f\sq\dom F,$
$F(u_\xi)=v_\xi$ for all $\xi,$ and still
$D=\dom F$ and $Z=\ran F$ consist of powers of $2$
and $\fras{f(d)}{d}\le\fras{f(d')}{d'}$ for all
$d<d'$ in~$D.$

Let $h=\cin D=\cin Z$ and
$D=\ans{d_1,d_2,\dots,d_h},$ 
$Z=\ans{z_1,z_2,\dots,z_h},$ 
in the increasing order of $\adN$ in $\dI.$
Then $z_\nu=F(d_\nu)$ for all $\nu=1,\dots,h.$
As all $d_\nu,z_\nu$ are powers of $2,$ the fractions
$j_\nu=\fras{d_\nup}{d_\nu}$ and
$k_\nu=\fras{z_\nup}{z_\nu}$
belong to $\adN$ and $j_\nu\le k_\nu$ by the above.
Note also that $d_1=u_0=1$.

Any number $x\in\adN$ admits, in $\dI,$ a unique
representation in the form
$x=\sum_{\nu=1}^h \al_\nu d_\nu,$ where
$\al_\nu\in\adN$ and 
$0\le \al_\nu<j_\nu$ for all $\nu=1,\dots,h-1$
(but $\al_h$ is not restricted, of course).
The first idea that comes to mind is to try
$\sg(x)=\sum_{\nu=1}^h \al_\nu z_\nu$ as a reduction of
$\meu$ to $\mev.$
However this does not work.
Indeed let $x=\sum_{\nu=1}^h d_\nu$ and
$x'=\sum_{\nu=1}^{h-1} (j_\nu-1) d_\nu,$  so that
$x-x'=1$ but $|\sg(x)-\sg(x')|$ can be very big
in the case when, say, $k_\nu>j_\nu$ for all $\nu.$
However there is a useful modification.

Suppose that $x=\sum_{\nu=1}^h \al_\nu d_\nu\in\adN,$
and $0\le \al_\nu<j_\nu$ for $\nu=1,\dots,h-1,$ as
above.
Say that $x$ is \rit{type-1}
if there exist indices $1\le \nu'<\nu''\le h-1$ such that
$d_{\nu'}\in U,$ $d_{\nu''}\nin U,$ and $a_\nu=j_\nu-1$
for all $\nu$ such that $\nu'\le\nu\le\nu''.$
Then take
the largest $\nu''$ and the least $\nu'$ such that the
pair $\nu',\nu''$ has this property, and put
$\bar \al_\nu=a_\nu$ for all $\nu<\nu'$ and $\nu>\nu'',$
$\bar \al_\nu=0$ for $\nu'\le\nu\le\nu'',$ and
$\bar \al_{\nu''+1}=\al_{\nu''+1}+1,$ and define
$\bar x=\sum_{\nu=1}^h \bar \al_\nu d_\nu.$
Otherwise ($x$ is \rit{type-2}) put $\bar x=x$.
Easily $\bar x-x=\al_{\nu'}\in U$ in the type-1 case.

\rit{Prove that the map\/ $\rho(x)=\sg(\bar x)$
is a reduction of\/ $\meu$ to\/ $\mev$},
that is,
${|x-x'|\in U}\eqv {|\sg(\bar x)-\sg(\bar y)|}\in V$
holds for all $x,x'\in\adN.$

Assume that
$x=\sum_{\nu=1}^h \al_\nu d_\nu$ and
$y=\sum_{\nu=1}^h \ga_\nu d_\nu,$
where $\al_\nu,\ga_\nu<j_\nu,$
and $|x-y|\in U,$ hence $|x-y|<u_\xi=d_{\nu}$
for some $\xi<\vt\zt\nu<h.$
Let $x<y.$ 
Assume \noo\ that $x,y$ are of type-2.
(Otherwise change $x,y$ to $\bar x,\bar y$.)
There exist infinitely (but \ddi finitely) many
indices $\nu'>\nu$ such that
$\al_{\nu'}\ne j_{\nu'}-1.$
In this case  
$\al_{\nu'}=\ga_{\nu'}$ for all $\nu'\ge\nu$
by the assumption $|x-y|<d_\nu.$   
Thus $|\sg(x)-\sg(y)|\in V$
(since $j_\nu\le k_\nu$ for all $\nu$), as required.

Now suppose that $x<y$ are as above, in particular,
of type-2, but $|x-y|\nin U,$
hence $|x-y|>u_\xi$ for all $\xi<\vt.$
Then 
$D'=\ens{d_\nu\in D}{\al_\nu\ne \ga_\nu}$
is an internal set, hence it has the largest element,
say $d_{\nu''}=\tmax D'.$ 
Note that $d_{\nu''}\nin U.$
(Use the assumption $|x-y|>u$ for all $u\in U$.)
We have $\al_{\nu''}<\ga_{\nu''}$ (as $x<y$).
Then the only opportunity for
$|\sg(x)-\sg(y)|$ to belong to $V$ is obviously
the existence of an index $\nu'<\nu''$ such that
$z_{\nu'}\in V$ and $\ga_\nu=0,$
$\al_\nu=j_\nu-1=k_\nu-1$ for all $\nu$ between
$\nu'$ and $\nu''.$
But this contradicts the assumption that $x$ is of
type-2.
Thus $|\sg(x)-\sg(y)|\nin V,$ as required.\vom

\ref{T2i} 
The proof of this item follows the same line as the
proof of \ref{T2f}, but with appropriate changes,
of course.
It will appear elsewhere.\vom

\ref{T4} 
Suppose that $U,V$ are resp.\ \ssi\ cofinal,
\ssi\ coinitial additive fast cuts. 
Prove that $\cef\meu\not\le\cef\mev;$ the proof of
$\cef\mev\not\le\cef\meu$ is similar.
Choose an increasing sequence $\sis{u_\xi}{\xi<\vt}$
and a decreasing sequence $\sis{v_\eta}{\eta<\tau}$ 
resp.\ cofinal in $U$ and coinitial in $\adN\bez V;$
$\vt=\cof U$ and $\tau=\coi V$ being
infinite regular cardinals in $\dV.$

Suppose on the contrary that  $R\sq \adN\ti\adN$
is an invariant pre-injection of $\meu$ to $\mev,$ 
that is,
${|x-x'|\in U}\eqv {|y-y'|}\in V$
for any pairs $\ang{x,y}$ and $\ang{x',y'}$
in $R,$ and $\dom R=\adN.$ 
Then $R=\bigcup_{a\in A}\bigcap_{b\in B}R_{ab},$ 
where $A,B\in\dV$ and $R_{ab}$ are internal sets.
Arguing as above in the proof of \ref{T2f}
(parts 1,2),
we obtain by \Sat\ 
for any fixed $a\in A$:
\envur{
\label{er1}
\bay[b]{l}
\hspace*{-2ex}
\sus\text{ finite }F\sq B\;\:
\sus \xi<\vt\;\:\sus\eta<\tau\quad
\kaz x,x',y,y'\in\adN:\\[0.4\dxii]
\qquad{\otn x{R_{aF}}y}\,\land\, 
{\otn{x'}{R_{aF}}{y'}}
\,\land\,{|y-y'|<v_\eta}\,\limp\,
{|x-x'|<u_{\xi}}\,,
\eay
}
where $R_{aF}=\bigcap_{b\in F}R_{ab},$ 
and, in  the opposite direction, 
\envur{
\label{er2}
\bay[b]{l}
\hspace*{-2ex}
\kaz\xi<\vt\;\:\kaz \eta<\tau\quad
\sus\text{ finite }F'\sq B\;\quad
\kaz x,x',y,y'\in\adN:\\[0.4\dxii]
\qquad{\otn x{R_{aF'}}y}\,\land\, 
{\otn{x'}{R_{aF'}}{y'}}
\,\land\,{|x-x'|<u_\xi}\,\limp\,
{|y-y'|<v_{\eta}}\,.
\eay
}

Let $a\in A.$
Take $\xi,\eta,F$ as in \eqref{er1}.
Take then $F'$ as in \eqref{er2} for $\xip$ and $\eta.$
We may assume that $F\sq F'.$
Then for all pairs 
$\ang{x,y}\zt\ang{x',y'}$ in the set $D(a)=\dom {R(a)},$
where $R(a)=R_{aF'},$ we have
$|x-x'|<u_{\xip}\limp{|x-x'|<v_{\xi}}.$ 
Assuming \noo\ that $\fras{u_{\xip}}{u_\xi}$
is infinitely large for all $\xi,$ we conclude that  
each $D(a)$ is an internal scattered set in the sense
of \ref{scar}, and so on towards contradiction
as above.\vom

\qeG{Thm~\ref T}

\seci{On small and large effective sets}
\label{fm}

Here we prove a \lap{small--large} type theorem 
related to $\DD$ quotients.
The notions of smallness and largeness will be
connected with a cut $U\sq\aCard,$ as in
Corollary~\ref{gc}.
By necessity there also will be a gap between the
largeness and smallness, but we don't know  
whether its size can be reduced.


Recall that a cut (initial segment) $U\sq\aCard$ is
called \rit{exponential}
iff ${\ka\in U}\imp{2^\ka\in U},$ or, equivalently,
$U=2^U$ holds, where
$2^U=\ens{\vt\in\aCard}{\sus \ka\in U\,(\vt\le2^\ka)}.$\,
($2^\ka$ is understood as the cardinal exponentiation
in $\dI$.)\,

We write $\la\ge2^U$ to mean $\la\ge2^\ka$ for all
$\ka\in U.$

\bte
\label{Ts}
Suppose that\/ $\qE$ is a\/ $\DD$ \eqr\ on an internal  
set\/ $\hh$ and\/ $U\sq\aCard$ is a\/ $\DD$ 
cut such that\/ $\dN\sq U.$
Then at least one of the following conditions holds\,$:$
\ben
\tenu{{\krm(\Alph{enumi})}}
\atlh
\itla{ls1}
for any\/ $\la\in\aCard$ with\/ $\la\ge 2^U$ and
any\/ $m\in\adN\bez\dN$
there is an internal
map\/ $\rho$ defined on\/ $\hh$ such that\/
$\cin{\ran\rho}\le \la^m$
{\krm(=$\la$ whenever $\la\nin\adN$)}
and\/ ${\rho(x)=\rho(y)}\imp{x\qE y}$
for all\/ $x,y\in\hh\,;$  

\itla{ls2}
there exists an internal set\/ $Y\sq\hh$ of pairwise\/
\dd\qE inequivalent elements such that\/ $\cin Y\nin U$.
\een
If\/ $U$ is an exponential non-internal cut then\/
\ref{ls1} and\/ \ref{ls2} are incompatible
even in the case when\/
$\DD$ maps\/ $\rho$ are allowed in\/ \ref{ls1}.
\ete

In terms of effective cardinals \ref{ls2} means
$\ka\le \cef{\hh/{\qE}}$
(and even by means of an internal reduction)
for some $\ka=\cin Y\in\adN\bez U,$ that
is a restriction of the cardinality of the quotient
$\hh/{\qE}$ from below.
Accordingly \ref{ls1} means that for all $\la\ge2^U$
and $m\in\adN\bez\dN$ and any internal $Z$ with
$\cin Z=\la^m$ there is an \eqr\ $\qF$ on $Z$
(in terms of \ref{ls1}, $\rho(x)\qF\rho(y)$ iff
$x\qE y$)
such that $\cef{\hh/{\qE}}\le\cef{Z/{\qF}}$
(still by means of an internal reduction),
a restriction of the cardinality of $\hh/{\qE}$
from above. 

Some theorems of this form are
known from descriptive set theory, for instance
Silver's theorem on $\fp11$ \eqr s in \ssyl{sil},
in which \lap{small} means at most countably many
equivalence classes while \lap{large} means that
there exists a pairwise \dd\qE inequivalent perfect
set.

Note that the implication
${\rho(x)=\rho(y)}\imp{x\qE y}$ in \ref{ls1}
cannot be replaced by the equivalence
${\rho(x)=\rho(y)}\eqv{x\qE y}$:
indeed the latter would imply the $\DD$ smoothness
of $\qE,$ which, generally speaking, is not the case
even for equivalence relations of the form $\qeu$
by Theorem~\ref{T}.\tsp

\noi
{\bf Proof} (Theorem~\ref{Ts}).
\rit{Case 1}: 
$U$ is standard size cofinal, including internal cuts. 
In this case we prove an even stronger result, namely 
the disjunction $\text{\ref{Ts1}}\lor\text{\ref{ls2}},$
where 
\ben
\tenu{{\mtho$(\text{\krm\Alph{enumi}}')$}}
\itla{Ts1}
there exist a set\/ $D\in\dV,$ and for each\/
$d\in D$ an internal set\/ $R_d$ and an internal map\/  
$f_d:\hh\to R_d$ such that\/ $\cin{R_d}\in 2^U$ and\/  
${f(x)=f(y)}\imp{x\qE y}$ for all\/ $x,y\in\hh,$ 
where\/ $f(x)=\sis{f_d(x)}{d\in D}\,.$
\een

We first show that \ref{Ts1} implies \ref{ls1}.
Suppose that
$\la\ge 2^U\zt m\in\adN\bez\dN.$
Recall that the map $d\mto\xd$ is an injection
$D\to \sD.$
Its image $D'=\ens{\xd}{d\in D}\sq\sD$ is a
set of standard size  together with $D.$
By \ref{lss1}(iii), $D'$ can be covered by an
internal set $S\sq \sD$ such that $\cin S\le m.$
The \Exn\ principle (1.3.13 in \ssyl{book})
yields an internal function $F$ defined on $S\ti\hh$
so that $F(\xd,x)=f_d(x)$ for all
$d\in D\zt x\in\hh.$
By the same reasons there is an internal map $r$
defined on $S$ so that $r(\xd)=R_d$ for all $d\in D$.
We can assume that for any $s\in S,$ $r(s)$
is an internal set with $\cin {r(s)}<\la,$
and $F(s,x)\in r(s)$ for all $x\in\hh.$
(Otherwise redefine $r$ and $F$ by
$r(s)=\ans0$ and $F(s,x)=0$ for all \lap{bad} $s$
--- but none of $s=\xd\zt d\in D,$ is \lap{bad}
in the assumptions of \ref{Ts1}.)
Put $\rho(x)(s)=F(s,x)$ for $x\in\hh\zt s\in S.$
 
We begin the proof of
$\text{\ref{Ts1}}\lor\text{\ref{ls2}}.$ 
By definition  
${\qE}=\bigcup_{a\in A}\bigcap_{b\in B}E^a_b,$
where all sets $E^a_b\sq\hh\ti\hh$ are internal while
$A,B\in\dV.$
We may \noo\ assume that every set $E^a_b$
is symmetric (similarly to $\qE$ itself),
that is, $E^a_b={(E^a_b)}\obr,$ where
$E\obr=\ens{\ang{y,x}}{\otn x{E}y}:$ 
indeed
\dm
\TS
{\qE}={\qE}\cap{{\qE}\obr}=
\bigcup_{a\in A}\;\bigcap_{b,b'\in B}\;
E^a_b\cup {(E^a_{b'})}\obr=
\bigcup_{a\in A}\;\bigcap_{b,b'\in B}\;C^a_{bb'}\,, 
\dm
where the sets
$C^a_{bb'}=(E^a_b\cup {(E^a_{b'})}\obr)\cap 
(E^a_{b'}\cup {(E^a_b)}\obr)$
are symmetric.
(We write $\otn x E y$ for 
$\ang{x,y}\in E$ whenever $E$ is a binary relation.)

It follows from the transitivity of $\qE$ that  
for any $x\zi y\in\hh$  
\dm
\sus a\in A\;\sus z\in\hh\;\kaz b\in B\;
({\otn x{E^a_b}z} \land {\otn y{E^a_b}z})
\limp {\otn x{\qE}y} \,. 
\dm 
The axiom of \Sat\ transforms this to
\dm
\sus a\in A\;\kaz B'\in \pfin(B)\;\sus z\in\hh\; 
({\otn x{E^a_{B'}}z} \land {\otn y{E^a_{B'}}z})
\limp {\otn x{\qE}y} \,, 
\dm
where $E^a_{B'}=\bigcap_{b\in B'}E^a_{b}.$ 
As the two leftmost quantifiers are
restricted to the sets  
$A$ and $\pfin(B)$ in $\dV,$
the last formula is equivalent to
\envur{
\label{zz1}
\kaz \vpi\in \Phi\;\sus a\in A\;\sus z\in\hh\; 
({\otn x{E^a_{\vpi(a)}}z} \land {\otn y{E^a_{\vpi(a)}}z})
\limp {\otn x{\qE}y} \,, 
}
where $\Phi\in\dV$ is the set of all functions
$\vpi:A\to\pfin(B)$.

As $U$ is standard size cofinal, there is an increasing
sequence $\sis{\nu_\xi}{\xi<\vt}$ of elements 
$\nu_\xi\in U,$ cofinal in $U,$  with $\vt$ 
being an infinite cardinal in $\Card,$ or simply $U$ 
is internal, $\vt=1=\ans0,$ and $\nu_0$ is the least 
element in $\aCard\bez U.$

Suppose that \ref{ls2} of the theorem fails, \ie\
there is no pairwise \dd\qE inequivalent sets $Y$ with
$\cin Y\nin U.$
More formally,
\dm
\kaz Y\in P\:
\big(
\kaz\xi<\vt\:(\cin Y\ge\nu_\xi)\imp \sus x\ne y \in Y \: 
\sus a\in A \: 
\kaz b\in B \: (\otn x{E^a_{b}}y)
\big)\,,
\dm
where $P=\pvn{H}=\ens{Y\sq\hh}{Y\text{ is internal}}.$ 
\Sat\ converts the expression to the right of $\imp$ to 
\dm
\sus a\in A \: \kaz B'\in\pfin(B)\:\sus x\ne y \in Y \: 
(\otn x{E^a_{B'}}y)\,,
\dm
and then to  
$ 
\kaz \vpi\in \Phi\;\sus a\in A\: 
\sus x\ne y \in Y \: 
(\otn x{E^a_{\vpi(a)}}y). 
$
We conclude that for any function $\vpi\in \Phi$  
\dm
\bay{l}
\kaz Y\in P\:
\big(
\kaz\xi<\vt\:(\cin Y\ge\nu_\xi)\imp \sus a\in A\: 
\sus x\ne y \in Y \: 
(\otn x{E^a_{\vpi(a)}}y) 
\big)\,.
\eay
\dm
\Sat\ yields an ordinal $\xi(\vpi)<\vt$ and a finite
set $A_\vpi\sq A$ such that
\envur{
\label{zz2}
\kaz Y\in P\:
\big(
\cin Y\ge\nu_{\xi(\vpi)}\imp \sus a\in A_\vpi\: 
\sus x\ne y \in Y \: 
(\otn x{E^a_{\vpi(a)}}y) 
\big)\,.
}
Let $Y_\vpi$ be any maximal
(internal) subset of $\hh$ such that
$\neg\;\otn x{E^a_{\vpi(a)}}y$ for all $a\in A_\vpi$  
and $x\ne y \in Y_\vpi.$ 
Then \eqref{zz2} implies
$\cin {Y_\vpi}<\nu_{\xi(\vpi)},$
while the properties of maximality of $Y_\vpi$
and symmetricity of $E^a_{b}$ 
imply 
\envur{
\label{pop3}
\kaz x\in\hh\;\sus y\in Y_\vpi\;\sus a\in A_\vpi\; 
(\otn x{E^a_{\vpi(a)}}y) \,.
}

Put 
$
\za_x(\vpi,a)=\ens{y\in Y_\vpi}
{\otn x{E^a_{\vpi(a)}}y} 
$
for $x\in \hh\zt\vpi\in\Phi\zt a\in A_\vpi.$
Thus $\za_x$ belongs to the set $Z$ of all functions 
$\za$ defined on the set
$D=\ens{\ang{\vpi,a}}{\vpi\in\Phi\land a\in A_\vpi}\in\dV$
and satisfying $\za_x(\vpi,a)\in R_\vpi=\pvn{Y_\vpi}.$
The sets $R_\vpi$ are internal and satisfy   
$\cin{R_\vpi}\in 2^U$ (because $\cin{Y_\vpi}\in U$).

We claim that $\za_x=\za_y$ implies $\otn x{\qE}y.$
It suffices, by \eqref{zz1}, to prove that for every
$\vpi\in\Phi$ there exist $a\in A\zt z\in \hh$
such that $\otn x{E^a_{\vpi(a)}}z$
and $\otn y{E^a_{\vpi(a)}}z.$ 
Note that $\za_x(\vpi,a)=\za_y(\vpi,a)\ne\pu$
for some $a\in A_\vpi$ by \eqref{pop3}. 
Take any $z\in\za_x(\vpi,a).$   
Then $z\in Y_\vpi,$ thus $\ang{x,z}$ and
$\ang{y,z}$ belong to $E^a_{\vpi(a)},$ as required.

To accomplish the proof of \ref{Ts1} in the
assumption $\neg\:\ref{ls2},$ we put
$f_{d}(x)=\za_x(\vpi,a)$ and $R_d=R_\vpi$ for all  
$x\in\hh$ and $d=\ang{\vpi,a}\in D.$\vom

\rit{Case 2}:
$U$ is standard size coinitial, but non-internal.
Suppose that \ref{ls2} fails, and consider any
$m\in\adN\bez\dN$ and $\la\ge 2^U.$ 
Then \ref{ls2} fails also fot the internal, hence,
\ssi\ cofinal, cut
$U'=\ens{\ka\in\aCard}{2^\ka\le\la}:$
indeed, $U\sq U'$ by the choice of $\la.$
Therefore \ref{ls1} holds for $U'.$
Thus there is an internal map $\rho\zt\dom\rho=\hh,$
such that $\cin{\ran\rho}\le \la^m$
and\/ ${\rho(x)=\rho(y)}\imp{x\qE y}.$\vom

\rit{Incompatibility}. 
Assume that $Y\sq\hh$ witnesses \ref{ls2},
in particular, $\ka=\cin Y\nin U=2^U.$
Then $U'=\ens{\la\in\aCard}{2^\la<\ka}$
is an internal cut with $U\sq U'.$
Thus $U\sneq U'$ since $U$ is non-internal.
Therefore there is $\vt\nin U$ such that
$2^\vt<\ka.$
Applying this trick once again, we find
$\vt\nin U$ with $2^{2^\vt}<\ka.$ 
Suppose on the contrary that $\rho$ witnesses
\ref{ls1} for $\la=2^\vt$ and some
$m\in\adN\bez\dN\zt m<\vt.$ 
Then $\rho\res Y$ is an internal injection of 
$Y$ into an internal set $Z=\imj\rho Y$
satisfying $\cin Z\le 2^{\vt m}.$
But this contradicts Theorem~\ref{ie},
since by definition
$2^{\vt m}\cdot n<2^{\vt\cdot\vt}<{2^{2^\vt}}<\ka=\cin Y$
for any $n\in\dN.$\vom

\qeG{Thm~\ref{Ts}}\tsp

The case $U=\dN$ deserves special attention.
Since $\dN$ is a \ssi\ cofinal cut, a stronger
dichotomy holds:
$\text{\ref{Ts1}}\lor\text{\ref{ls2}}.$ 
Clearly \ref{ls2} claims 
the existence of an infinite internal set 
of pairwise \dd\qE inequivalent elements in this case.
On the other hand, the sets $R_d$ in \ref{Ts1}
are finite, hence $P=\prod_{d\in D}R_d$ is a set of
standard size, and so is any quotient of the form
$P/{\qF},$ where $\qF$ is an \eqr\ on $P.$ 
Thus \ref{Ts1} implies that
$\hh/{\qE}$ itself is a set of standard size. 
Such a dichotomy
(\ie\ standard size of $\hh/{\qE}$ or
an infinite internal pairwise inequivalent set)
is contained in Theorem 1.4.11 in \ssyl{book}. 
Similar dichotomies appeared in \ssyl{mon} for
countably determined \eqr s.
P.\,Zlato\v{s} informed us that a close  
result for $U=\dN$ was earlier obtained by Vencovsk\'a 
(unpublished)
in the frameworks of \AST.

\seci{Nonstandard version of the finite Ramsey theorem}
\label{part}

The following corollary of Theorem~\ref{Ts} is a
Ramsey--like result.
Recall that $[A]^n=\ens{X\sq A}{\card X=n}.$
By a \rit{partition} of $[A]^n$ we understand any
\eqr\ $\qE$ on $[A]^n,$ and a \rit{homogeneous set}
for $\qE$ is any $H\sq A$ such that the sets
$X\in[H]^n$ are pairwise \dd\qE equivalent.

The finite Ramsey theorem claims (in \ZFC) that
\ben
\fenu
\addtocounter{enuf}{-1}
\itla{frt}
\rit{for any natural numbers\/
$\ell,n,s$ there is\/ $k\in\dN$ such that\/
$k\to(\ell)_s^n$.}
\een
Here $k\to(\ell)_s^n$ means that
for any partition of $[k]^n$
into \dd s{}many parts there is an \dd\ell element
homogeneous set $H\sq k.$
We refer to \ssyl r, and also to  
3.3.7 in \ssyl{cc}, \S\,6 in \ssyl{kun}, or \ssyl{ehmr} 
for a modern proof, details and related results. 

Let $K(\ell,s,n)$ denote the least $k$ satisfying 
$k\to(\ell)_s^n.$
It is known that 
$K(\ell,s,n)$ is rapidly increasing as a function of 
$\ell$ for any fixed $n,s,$ see \ssyl{ehmr}.
But of course $K$ is a recursive function. 

It is an easy nonstandard corollary of \ref{frt}
that $\kpa\arin(\ell)_s^n$ for all $n,s,\ell\in\dN$ and
$\kpa\in\adN\bez\dN$ where {\tt int} over the arrow means
that the partition and the homogeneous set
are assumed to be internal.
A nicer nonstandard version, also well-known, is
\rit{$\kpa\arin(\iy)_s^n$ for any\/ $n,s\in\dN$ and\/
$\kpa\in\adN\bez\dN$}, that is,  
any internal partition\/ $[\kpa]^n$ into\/ $s$ parts 
admits an infinite internal homogeneous set.
By the way, its quantifier structure is simpler than
that~of~\ref{frt}:\pagebreak[0]
\dm
\kaz \kpa,\ell,n,s\;\kaz\,\text{partition}\;
\sus A\;\kaz u\zi v\in [A]^n.
\dm 

The following theorem contains a much more general
claim.
In \HST, define a function $K$ in $\dV$ as above.
Then $\sK$ is a standard function $\adN^3\to\adN$ having 
in the internal universe $\dI$ the same
properties as $K$ in~$\dV.$

\bte
\label{mlzl}
Suppose that\/ 
$U\sneq\adN$ is a\/ $\DD$ cut with\/ $\dN\sq U,$
closed under\/ $\sK$ and exponential,
$n\in\dN\zt \kpa\in\adN\bez U,$
and\/ $\qE$ is a\/ $\DD$ \eqr\ on\/ $[\kpa]^n.$
If there is no internal pairwise \dd\qE inequivalent
sets\/ $Y\sq[\kpa]^n$ satisfying\/ $\cin Y\nin U,$
then the partition\/ $\qE$ admits an internal
homogeneous set\/ $A\sq \kpa$ such that\/
$\cin A\nin U.$ 
\ete

A similar result was obtained in \ssyl{mz'} in the case
$U=\dN$ for countably determined \eqr s.
See Theorem~2.8 in \ssyl{kkml} for a somewhat weaker
result in the case when $t$ in the proof of \ref{mlzl}
is predefined.

\bpF
Define, in $\dV,$ $f(s)=K(s,s,n)$  for each $s\in\dN.$
Then $f:\dN\to\dN$ and $s\le f(s)\zd\kaz s.$
The map $\stf$ has the same properties with respect to 
$\adN.$
As $U$ is  \dd\sK closed and exponential, there exist
$s,\vt\in\adN\bez U$ and $m\in\adN\bez\dN$
such that $\stf(s)=\sK(s,s,n)\le\kpa$
and $2^{\vt m}\le s.$  

In our assumptions, \ref{ls2} of Theorem~\ref{Ts} fails,
hence \ref{ls1} holds, that is, there exists an internal
map $\rho$ defined on\/ $[\kpa]^n$ such that\/
$\cin{\ran\rho}\le 2^{\vt m}\le s$
and\/ ${\rho(u)=\rho(u)}\imp{u\qE u}$
for all\/ $u,v\in[\kpa]^n.$  
On the other hand, we have $\kpa\to(s)^n_s$ by the
choice of $s,$ therefore the partition of 
$[\kpa]^n$ induced by $\rho$ has an internal homogeneous
set $A$ such that $\cin A=s\nin U.$
Thus $\rho(u)=\rho(v),$ and hence $u\qE v,$
for all $u,v\in[A]^n.$
\epF

\def\refname{{\large\bf References}}

\end{document}